\documentclass[11pt,reqno]{amsart}

\usepackage{color,xcolor}
\usepackage{enumitem}
\usepackage{soul}
\usepackage{array}
\usepackage{cite}
\usepackage{amsmath, amssymb, amsthm}
\usepackage{multirow, multicol}
\usepackage{caption}
\usepackage{scalerel}
\usepackage{algorithm}
\usepackage{algpseudocode}
\usepackage{url}
\usepackage{longtable}
\usepackage{relsize}
\usepackage{kotex}
\usepackage[colorinlistoftodos]{todonotes}

\newcommand{\commentout}[1]{}

\def\a{{\alpha}}

\def\d{{\delta}}

\def\k{{\kappa}}

\def\l{{\lambda}}
\def\z{{\zeta}}

\def\ff{{\mathcal F}}

\def\F{{\mathbb F}}

\def\R{{\mathbb R}}
\def\Z{{\mathbb Z}}
\def\C{{\mathbb C}}

\def\HE{{\mathcal {HE}}}

\DeclareMathOperator\re{Re}

\DeclareMathOperator\av{av}

\newtheorem{thm}{Theorem}[section]

\newtheorem{lem}[thm]{Lemma}

\theoremstyle{definition}
\newtheorem{defn}[thm]{Definition}
\newtheorem{rem}[thm]{Remark}
\newtheorem{ex}[thm]{Example}

\newtheorem{innercustomdef}{Definition}
\newenvironment{customdef}[1]
{\renewcommand\theinnercustomdef{#1}\innercustomdef}
{\endinnercustomdef}

\newtheorem{innercustomthm}{Theorem}
\newenvironment{customthm}[1]
{\renewcommand\theinnercustomthm{#1}\innercustomthm}
{\endinnercustomthm}

\newtheorem{innercustomcor}{Corollary}
\newenvironment{customcor}[1]
{\renewcommand\theinnercustomcor{#1}\innercustomcor}
{\endinnercustomcor}

\begin{document}

\title[The $L$-polynomial of hyperelliptic function fields]
{The $L$-polynomial of hyperelliptic function fields\\ and its applications}

\author{Peter J. Cho and Jinjoo Yoo$^*$}

\address{Department of Mathematical Sciences, Ulsan National Institute of Science and Technology, 50, UNIST-gil, Ulsan, 44919, South Korea,
\ E-mail : petercho@unist.ac.kr}

\address{Center for Quantum Structures in Modules and Spaces, Seoul National University, 1, Gwanak-ro, Gwanak-gu, Seoul, 08826, South Korea,
\ E-mail : jinjooyoo121@gmail.com}

\date{}


\thanks{
{\tiny {${}^{*}$ Corresponding author.\\
E-mail address: petercho@unist.ac.kr (Peter J. Cho), jinjooyoo121@gmail.com (J. Yoo)\\
The first author is supported by the National Research Foundation of Korea (NRF) grant funded by the Korea government (MSIT) (No. RS-2022-NR069491 and No. RS-2025-02262988). The second author is supported by the NRF grant funded by the Korea government (No. 2020R1A5A1016126 and RS-2022-NR069491).}}}

\begin{abstract}
Let $\ell$ be an odd prime, $q$ an odd prime power such that $q \not\equiv 0 \pmod \ell$, and $m$ the order of $q$ in $\F_\ell^\times$. We propose an explicit $L$-polynomial of hyperelliptic function field $K:=\F_q(T, \sqrt[\ell]{T^2+aT+b})$ with $a, b \in \F_q$ and $a^2-4b \ne 0$. Using our formula, we obtain the explicit closed formula for the class number of $K$, where $m$ is even or $m=\frac{\ell-1}{2}$.
As an application, we compute the average class numbers for hyperelliptic function fields with genus up to $3$.
\end{abstract}

\maketitle


\section{Introduction}

The Dedekind zeta function $\zeta_F(s)$ of a number field $F$ is a fundamental tool for studying the arithmetic properties of number fields, including class numbers and the distribution of prime ideals. For the global function field case, the zeta function is defined as follows.
$$\z_K(s) = \sum_{A\ge 0} NA^{-s} \quad (\re(s)>1),$$
where $K$ is an algebraic function field of one variable over a finite field $\F_q$, $A \in \mathcal{D}_K$, and $NA := q^{\deg A}$ with divisor group $\mathcal{D}_K$ of $K$. We note that $\mathcal{D}_K$ is the (additively written) free abelian group which is generated by the places of $K$.

It is well known that the zeta function $\z_K(s)$ of a global function field $K$ can be expressed as {\cite[Theorem 5.9]{R}}
\begin{equation}\label{zeta}
\z_K(s) = \frac{L_K(q^{-s})}{(1-q^{-s})(1-q^{1-s})};
\end{equation}
this holds for all $s$ such that $\re(s)>1$ and the right-hand side provides an analytic continuation of $\z_K(s)$ to all of $\C$. The zeta function $\z_K(s)$ has simple poles at $s = 0$ and $s = 1$. The coefficients of $L_K(q^{-s})$ are integer and the degree of $L_K(q^{-s})$ is $2g$, where $g$ is the genus of $K$.

For convenience, we let $u:=q^{-s}$. Then by \eqref{zeta}, the zeta function $\z_K(s)$ of $K$ can be written as in the form of
$$\z_K(s) = \frac{L_K(u)}{(1-u)(1-qu)}, \text{ where } L_K(u) \in \Z[u].$$
We call the numerator of the zeta function $L_K(u) = \sum_{i=0}^{2g}c_iu^i \in \Z[u]$ by the \textit{$L$-polynomial of $K$}. The $L$-polynomial of a global function field $K$ encodes fundamental arithmetic and geometric information about $K$. In particular, it determines the number of points of curves over finite field, determines the eigenvalues of Frobenius, and reflects the structure of the Jacobian of the associated curve. Consequently, understanding the $L$-polynomial of $K$ is an important problem in the study of global function fields.

One of the applications for the $L$-polynomial of $K$ is the \textit{class number}. Computation of the class number for a number field $F$ of large degree is a hard task because we need to find an integral basis and fundamental units for $F$, and also determine the residue of the Dedekind zeta function $\z_F(s)$ at $s=1$. For the case of the global function field $K$, the class number $h_K$ of $K$ is given by $L_K(1)$. Letting $L_K(u) = \sum_{i=0}^{2g}c_iu^i$, the class number \mbox{$h_K$ of $K$ is}
\begin{equation*}
h_K = (q^g+1) + (q^{g-1}+1)c_1 + (q^{g-2}+1)c_2 + \cdots + (q+1)c_{g-1}+ c_g;
\end{equation*}
we use the relation between coefficients $c_i$'s (see Lemma \ref{b}). We note that the class number $h_K$ tends to increase as the genus $g$ of $K$ grows. The main obstacle for computing $L_K(1)$ value is the determination of the number $\sharp \mathcal{C}(\F_{q^t})$ of the rational points of the curve $\mathcal{C}$ over $\F_{q^t}$ with $t \le g$, where $g$ is the genus of $\mathcal{C}$ and $\mathcal{C}$ is the curve which corresponds to the field $K$. Even though the classifications of function fields with class numbers less than equal to 3 are done in \cite{L1, L2, LM, MS, P, SS}, for the reason mentioned above, it is very challenging to compute the class number of a global function fields with large genus.

In this paper, we consider \textit{(hyper)elliptic function fields}. Let $\ell$ be an odd prime and $q$ an odd prime power with $q \not \equiv 0 \pmod \ell$. Let the set of (hyper)elliptic function fields be as follows:
$$\HE_\ell(\F_{q}):=\{k(\sqrt[\ell]{T^2+aT+b}) \mid a, b \in \F_q,~ a^2-4b \ne 0\} \text{ with $k:=\F_q(T)$.}$$
For the case $a^2-4b=0$, the genus of $k(\sqrt[\ell]{T^2+aT+b})$ is zero; this implies that the $L$-polynomial of $k(\sqrt[\ell]{T^2+aT+b})$ is one.

We introduce the definition of Jacobi sum, which plays an important role in our main result. We note that a multiplicative character $\l$ of $\F_q$ of order $\ell$ can only exist when $q \equiv 1 \pmod \ell$.

\begin{customdef}{1}
Let $\l_1, \ldots, \l_n$ be multiplicative characters of $\F_q$ with $\l_i(0):=0$ $(1 \le i \le m)$. Then the sum
$$J(\l_1,\ldots, \l_n) = \sum_{c_1 + \cdots + c_m=1} \l_1(c_1)\cdots\l_n(c_n)$$
with the summation extended over all $n$-tuples $(c_1, \ldots, c_n)$ of elements of $\F_q$ satisfying $c_1 + \cdots + c_n = 1$ is called a \textit{Jacobi sum} in $\F_q$.
\end{customdef}

\begin{rem}
Let $J(\l_1, \ldots, \l_n) = \sum_{i=1}^\ell a_i\z^i$, where $\z = e^{\frac{2\pi \imath}{\ell}}$ and $a_i \in \Z$. Using the relation $\z^{\ell-1} = -(1+\z+\cdots + \z^{\ell-2})$, a Jacobi sum $J(\l_1, \ldots, \l_n)$ can be uniquely represented as $\sum_{i=0}^{\ell-2}\overline{a_i}\z^i$. However, in this paper, we use the Jacobi sum of the form $\sum_{i=1}^\ell a_i\z^i$: that is, the form without applying relation $\z^{\ell-1} = -(1+\z+\cdots + \z^{\ell-2})$.
\end{rem}

We now state our main results. Theorem \ref{thmA} determines the $L$-polynomial of the global function field $K \in \HE_\ell(\F_q)$. From this, Corollary \ref{corA} follows, giving the number of $\F_{q^t}$-points on the curve $HE_{a,b}$ associated with $K$. Finally, Theorem \ref{thmB} provides the explicit formulas for the $L$-polynomial and class number of $K$ with specific $m$.

\begin{customthm}{2}\label{thmA}
Let $\ell \ge 3$ be a prime, $q$ an odd prime power, and $K=\F_q(T, \sqrt[\ell]{T^2+aT+b})$ be a (hyper)elliptic function field over $\F_q$, where $T^2+aT+b\in \F_q[T]$ and $a^2-4b \ne 0$. Let $m$ be the order of $q$ in $\F_\ell^\times$. Then the $L$-polynomial of $K$ is
\begin{align*}
&L_K(u) = \sum_{i=0}^{\frac{\ell-1}{2}-1}c_i(u^i + q^{\frac{\ell-1}{2}-i}u^{\ell-1-i})+c_{\frac{\ell-1}{2}}u^{\frac{\ell-1}{2}} \quad \text{with $c_0=1$ and}\\[5pt]
&c_t = \frac{1}{t}\sum_{i=1}^tS_ic_{t-i} \text{ for $1 \le t \le \frac{\ell-1}{2}$, where}
\end{align*}
\begin{itemize}[leftmargin=0.2in]
\vskip 0.2cm
\item $\l_1$ (resp. $\l_2$) : multiplicative character of $\F_{q^m}$ of order $\ell$ (resp. $2$),
\vskip 0.1cm
\item $(-1)^{r-1}J(\l_1,\l_2)^r := \sum_{i=1}^\ell a_{r,i}\z^i$ with $\z:=e^{\frac{2\pi \imath}{\ell}}$,
\vskip 0.1cm
\item $\z^n=\l_1(-4^{-1}(a^2-4b))$, where $1\le n \le \ell$,
\vskip 0.1cm
\item $\ff_{r,n} := \ell a_{r,\overline{\ell-rn}}-\sum_{i=1}^\ell a_{r,i}$ $(1 \le n \le \ell)$ with $\overline{\ell-rn} \equiv \ell-rn \pmod \ell$ and $a_{r,0}:=a_{r,\ell}$,
\vskip 0.1cm
\item for $1\le t \le \frac{\ell-1}{2}$,
$$S_t = \left\{
\begin{array}{ll}
0 &\text{if $m \nmid t$ or $m > \frac{\ell-1}{2}$,}\\[3pt]
\ff_{r,n}& \text{if $t=mr$ and $a^2-4b =\square$ in $\F_{q^m}$,}\\[5pt]
(-1)^r\ff_{r,n}& \text{if $t=mr$ and $a^2-4b \ne \square$ in $\F_{q^m}$.}
\end{array}\right.$$
\end{itemize}
\end{customthm}

\vskip 0.3cm

\begin{customcor}{3}\label{corA}
Let all the notations be the same as in Theorem \ref{thmA} and
\begin{equation}\label{he}
HE_{a,b}: y^\ell = x^2+ax+b 
\end{equation}
be a (hyper)elliptic curve over $\F_q$, where $x^2+ax+b\in \F_q[x]$ and $a^2-4b \ne 0$. Let $m$ be the order of $q$ in $\F_\ell^\times$.

For a positive integer $t$, the number of $\F_{q^t}$-points $\#HE_{a,b}(\F_{q^t})$ of $HE_{a,b}$ is given as follows:
$$\left\{
\begin{array}{ll}
q^t+1 &\text{if $m \nmid t$ or $m > \frac{\ell-1}{2}$,}\\[3pt]
q^t+1+\ff_{r,n}& \text{if $t=mr$ and $a^2-4b =\square$ in $\F_{q^m}$,}\\[5pt]
q^t+1+(-1)^r\ff_{r,n}& \text{if $t=mr$ and $a^2-4b \ne \square$ in $\F_{q^m}$,}
\end{array}\right.$$
\end{customcor}

\vskip 0.3cm

\begin{customthm}{4}\label{thmB}
Let $\ell$ be an odd prime, $q$ an odd prime power, and $m$ the order of $q$ in $\F_\ell^\times$. Let $J(\l_1, \l_2) := \sum_{i=1}^\ell a_{1,i}\z^i$ with $\z = e^{\frac{2\pi \imath}{\ell}}$, where $\l_1$ (resp. $\l_2$) is a multiplicative character of $\F_{q^m}$ of order $\ell$ (resp. 2). The $L$-polynomial of $K=\F_q(T, \sqrt[\ell]{T^2+aT+b})$ with $a,b \in \F_q$ and $a^2-4b\ne 0$ is
\allowdisplaybreaks
\begin{equation*}
L_K(u) = \left\{\begin{array}{ll}
q^{\frac{\ell-1}{2}}u^{\ell-1}+1 &\text{if $m = \ell-1$}, \\[3pt]
\sum_{i=0}^{r-1}\binom{2r}{i}(q^{\frac{(2r-i)m}{2}}u^{(2r-i)m}+q^{\frac{im}{2}}u^{im})+\binom{2r}{r}q^{\frac{rm}{2}}u^{rm} &\text{if $m \mid \frac{\ell-1}{2}$ and $m$ is even},\\[5pt]
\sum_{i=0}^{r}\binom{2r+1}{i}(q^{\frac{(2r+1-i)m}{2}}u^{(2r+1-i)m}+q^{\frac{im}{2}}u^{im})&\text{if $m \nmid \frac{\ell-1}{2}$ and $m$ is even},\\[5pt]
q^{\frac{\ell-1}{2}}u^{\ell-1}+ 1 - \frac{\ell a_{1,\ell}+1}{m}u^{\frac{\ell-1}{2}} &\text{if $m = \frac{\ell-1}{2}$ is odd}\\[3pt]
&\hskip 0.5cm\text{and $a^2-4b = \square$ in $\F_{q^m}^\times$},\\[5pt]
q^{\frac{\ell-1}{2}}u^{\ell-1}+1+\frac{\ell a_{1,\ell}+1}{m}u^{\frac{\ell-1}{2}}&\text{if $m = \frac{\ell-1}{2}$ is odd}\\[3pt]
&\hskip 0.5cm \text{and $a^2-4b \ne \square$ in $\F_{q^m}^\times$},
\end{array}\right.
\end{equation*}
where $r:= \lfloor \frac{\ell-1}{2m}\rfloor$.

Furthermore, the class number of $K$ is
$$h_K=\left\{\begin{array}{ll}
(q^{\frac{m}{2}}+1)^{\frac{\ell-1}{m}} & \text{if $m$ is even,}\\[5pt]
q^m+1 - \frac{\ell a_{1,\ell} + 1}{m} & \text{if $m = \frac{\ell-1}{2}$ is odd and $a^2-4b = \square$ in $\F_{q^m}^\times$},\\[5pt]
q^m+1 + \frac{\ell a_{1,\ell} + 1}{m} & \text{if $m = \frac{\ell-1}{2}$ is odd and $a^2-4b \ne \square$ in $\F_{q^m}^\times$}
\end{array}\right.$$
with $a_{1, \ell} = 2N-\frac{q^m-1-\ell}{\ell}$, where $N$ is the number of quadratic residue elements in a set $T:=\{1-g^{\ell i} \mid 1 \le i \le \frac{q^m-1-\ell}{\ell}$\} and $g$ is a primitive root for $\F_{q^m}^\times$.
\end{customthm}

Using Theorem \ref{thmB}, we can determine the $L$-polynomial and the class number of a global function field with large genus. For example, let $K = \F_q(T, \sqrt[199]{T^2+aT+b})$ be a function field over $\F_q$ with $\ell=199$: that is, the genus of $K$ is 99. If we take $q=11$, then the order of 11 in $\F_{199}^\times$ is 22. By Theorem \ref{thmB}, the $L$-polynomial is
$$L_K(u)=\sum_{i=0}^4\binom{9}{i}(11^{11(9-i)}u^{22(9-i)}+11^{11i}u^{22i})$$
and the class number of $K$ is $(11^{11}+1)^9$ which has 104 digits. Letting $q=7$, we find that the order of 7 in $\F_{199}^\times$ is $99 = \frac{199-1}{2}$. Therefore, again by Theorem \ref{thmB}, the $L$-polynomial of $K$ is $7^{99}u^{198}+1\mp \frac{199a_{1,199}+1}{99}u^{99}$ and the class number of $K$ is $7^{99}+1 \mp\frac{199\cdot a_{1,199}+1}{99}$; the sign depends on whether $a^2-4b$ is a square or not in $\F_{7^{99}}^\times$. We note that Stoll \cite{S1, S2} investigated the arithmetic of the curves $C_A:y^2=x^\ell+A$ and \cite[Corollary 3.1]{S2} yields the case where $m$ is even.


\vskip 0.3cm


This paper is organized as follows. Section \ref{prelim} introduces some known facts related to function fields and finite fields. In Section \ref{sec 3}, we compute the explicit formulas for the trace of Frobenius of $HE_{a,b}$ over $\F_{q^t}$ (Theorem \ref{thm1}) and provide proofs of Theorem \ref{thmA}, Corollary \ref{corA}, and Theorem \ref{thmB}. Section 4 presents applications of Theorem \ref{thmA}, giving explicit formulas for the average class numbers of $\HE_5(\F_q)$ and $\HE_7(\F_{31})$ in Theorem \ref{thmC} and Theorem \ref{thmD}, respectively. Finally in Section 5, we provide Examples \ref{ex2} through \ref{class_num_ex4}, which illustrate Theorem \ref{thmA} and Corollary \ref{corA}. Example \ref{ex1} demonstrates the computation of the average class number on $\HE_5(\F_{31})$, related to Theorem \ref{thmC}.


\section{Preliminaries}\label{prelim}

As before, we let $K$ be an algebraic function field of one variable over a finite field $\F_q$. The following two lemmas state the well-known properties of the $L$-polynomial $L_K(u)$ of $K$.

\begin{lem}{\rm{\cite[Theorem 5.1.15]{St}}}\label{b}
Let $L_K(u)=\sum_{i=0}^{2g}c_iu^i \in \Z[u]$ be the $L$-polynomial of function field $K$ whose genus is $g$.
\begin{itemize}[leftmargin=0.3in]
\item [{\rm{(i)}}] $c_0 = 1$ and $c_{2g} = q^g$.
\item [{\rm{(ii)}}] $c_{2g-i} = q^{g-i}c_i$ for $0 \le i \le g$.
\item [{\rm{(iii)}}] $c_1 = N_1 - (q + 1)$, where $N_1$ is the number of places in $K$ whose degrees are one.
\item [{\rm{(iv)}}] $L_K(u)$ factors in $\C[u]$ in the form
\begin{equation}\label{L}
L_K(u) = \prod_{i=1}^{2g}(1-\a_iu).
\end{equation}
The complex numbers $\a_1, \ldots, \a_{2g}$ are algebraic integers, and they can be arranged
in such a way that $\a_i\a_{g+i} = q$ holds for $1 \le i \le g$.
\item [{\rm{(v)}}] If $L_{K_r}(u)$ denotes the $L$-polynomial of the constant
field extension $K_r:= K\F_{q^r}$, then
$$L_{K_r}(u) = \prod_{i=1}^{2g}(1-\a_i^ru),$$
where the $\a_i$ are given by \eqref{L}.
\end{itemize}
\end{lem}

\vskip 0.3cm

\begin{lem}{\rm{\cite[Corollary 5.1.17]{St}}}\label{c}
Let $L_K(u) = \sum_{i=0}^{2g}c_iu^i$ be the $L$-polynomial of $K/\F_q$, and $S_r := N_r - (q^r + 1)$, where $N_r$ is the number of rational places of $K_r:=K\F_{q^r}$. For $1 \le i \le g$, we have the following:
$$c_0=1 \text{ and } ic_i = S_ic_0 + S_{i-1}c_1 + \cdots + S_1c_{i-1}.$$
\end{lem}

We now give a definition of a \textit{hyperelliptic function field}.

\begin{defn}
A hyperelliptic function field over $\F_q$ is an algebraic function field $K/\F_q$ of genus $g \ge 2$ which contains a rational subfield $\F_q(T) \subseteq K$ with $[K : \F_q(T)] = 2$.
\end{defn}

A well-known property of a hyperelliptic function field is given as follows:

\begin{lem}{\rm{\cite[Lemma 6.2.2]{St}}}\label{hyper_prop}
\noindent
\begin{itemize}[leftmargin=0.2in]
\vspace{-0.25cm}
\item [{\rm{(i)}}] A function field $K/\F_q$ of genus $g \ge 2$ is hyperelliptic if
and only if there exists a divisor $A \in \mathcal{D}_K$ with $\deg A = 2$ and $\ell(A) \ge 2$.
\item [{\rm{(ii)}}] Every function field $K/\F_q$ of genus $2$ is hyperelliptic.
\end{itemize}
\end{lem}

Let $K=k(\sqrt[\ell]{Q(T)})$ be a cyclic extension over the rational function field $k:=\F_q(T)$, where $\ell$ is a prime number, $q$ is a prime power such that $q \not\equiv 0 \pmod \ell$, and $Q(T) \in \F_q[T]$. We denote $\infty$ by the \textit{infinite place of $k$} which corresponds to $(\frac{1}{T})$. If $q \equiv 1 \pmod \ell$ (respectively, $q \not \equiv 1 \pmod \ell$), then we call such extension as \textit{Kummer extension} (respectively, \textit{non-Kummer extension}).
A cyclic extension $K/k$ of degree $\ell$ is called an \textit{Artin-Schreier extension} if $q\equiv 0$ $\pmod \ell$. This extension is not defined by adjoining the $\ell$-th root of a polynomial $Q(T)$ to $k$, but in a different way. See \cite[p.256]{HL} for details.

\vskip 0.3cm

\begin{lem}\label{hyperelliptic_ff}
Let $K = k(\sqrt[\ell]{Q(T)})$ be a cyclic extension over $k=\F_q(T)$, where $\ell$ is an odd prime, $q$ is a prime power such that $\ell \nmid q$, and $Q(T)$ is a monic quadratic polynomial over $\F_q$. Then the ramified places of $k$ are the divisors of $Q(T)$ and $\infty$. Furthermore, the genus $g$ of $K$ is
$$g=\left\{
\begin{array}{ll}
0 & \text{if $Q(T) = (T-c)^2$ for some $c\in \F_q^\times$,}\\[3pt]
\frac{\ell-1}{2}& \text{otherwise.}
\end{array}
\right.$$
\end{lem}

\begin{proof}
This follows immediately from \cite[Proposition 3.7.3]{St} and \cite[Lemma 3.3]{LY}.
\end{proof}

We note that by Lemma \ref{hyperelliptic_ff}, the genus $g$ of $K=k(\sqrt[\ell]{Q(T)})$ with $\ell \ge 5$ is $g=\frac{\ell-1}{2} \ge 2$. Let $\mathfrak{p}_\infty$ be the place of $K$ which lies over the infinite place $\infty$ of $k$. By Lemma \ref{hyperelliptic_ff}, the degree of $\mathfrak{p}_\infty$ is one; this follows from the fact that $\infty$ of $k$ is ramified in $K$. Taking the divisor $A$ of $K$ as $A=2\mathfrak{p}_\infty$, we obtain that $\deg(A)=2$ and $\ell(A) \ge 2$. Thus, by Lemma \ref{hyper_prop}, $K$ is a hyperelliptic function field.

The rest of this section describes about Jacobi sum and the formula for estimating the number of solutions of a diagonal equation which are crucial for our main results.

\begin{defn}
A \textit{diagonal equation} over $\F_q$ is
$$a_1x_1^{k_1} + a_2x_2^{k_2} + \cdots + a_nx_n^{k_n}=b$$
with positive integers $k_1, \ldots, k_n$, $a_i \in \F_q^\times$ (for $1 \le i \le n$), and $b\in \F_q$.
\end{defn}

The number of solutions
$$N=N(a_1x_1^{k_1} + a_2x_2^{k_2} + \cdots + a_nx_n^{k_n}=b)$$ in $\F_q^n$ can be expressed in terms of Jacobi sums, where the definition is given as follows.

\vskip 0.4cm

\begin{lem}{\rm{\cite[Theorem 6.34]{LN}}}\label{number}
For $b\in \F_q^\times$, the number $N$ of solutions of the diagonal equation $a_1x_1^{k_1} + a_2x_2^{k_2} + \cdots + a_nx_n^{k_n}=b$ in $\F_q^n$ is given by
$$N=q^{n-1} + \sum_{j_1=1}^{d_1-1} \cdots \sum_{j_n=1}^{d_n-1}\lambda_1^{j_1}(ba_1^{-1})\cdots\lambda_n^{j_n}(ba_n^{-1})J(\lambda_1^{j_1}, \cdots, \lambda_n^{j_n}),$$
where $\lambda_i$ is a multiplicative character of $\F_q$ of order $d_i = \gcd (k_i,q-1)$.
\end{lem}

\begin{lem}{\rm{\cite[Theorem 5.26]{LN}}}\label{jacobi_lift}
Let $\l_1, \ldots, \l_k$ be multiplicative characters of $\F_q$, not all of which are trivial. Suppose $\l_1, \ldots, \l_k$ are lifted to characters $\l_1', \ldots, \l_k'$, respectively, of the finite extension field $E$ of $\F_q$ with $[E:\F_q]=s$. Then
\begin{equation*}
J(\l_1', \ldots, \l_k') = (-1)^{(s-1)(k-1)}J(\l_1, \ldots, \l_k)^s.
\end{equation*}
\end{lem}


\section{Determination for the $L$-polynomial of $k(\sqrt[\ell]{T^2+aT+b})$}\label{sec 3}

In this section, we state Theorem \ref{thm1} which gives the explicit formulas for the trace of Frobenius of $HE_{a,b}$ over $\F_{q^t}$, where $HE_{a,b}$ is defined in \eqref{he}. Combining Lemmas \ref{b}, \ref{c}, and Theorem \ref{thm1}, Theorem \ref{thmA} follows immediately. Corollary \ref{corA} then follows from Theorem \ref{thmA}. Finally, we provide the proof of Theorem \ref{thmB}. We recall that, by Lemma \ref{hyperelliptic_ff}, the genus $g$ of $K$ is $\left\{\begin{array}{ll}
0  & \text{if $a^2-4b =0$,}\\[3pt]
\frac{\ell-1}{2}  & \text{otherwise;}
\end{array} \right.$ and the $L$-polynomial of the function field with genus zero is always one.

\vskip 0.3cm

\begin{thm}\label{thm1}
Let $m$ be the order of $q$ in $\F_\ell^\times$, where $\ell$ is an odd prime and $q$ is an odd prime power. Let $HE_{a,b} : y^\ell = x^2+ax+b$ be a (hyper)elliptic curve which is defined in \eqref{he}. For a positive integer $t$, we denote the trace of Frobenius trace of $HE_{a,b}$ over $\F_{q^t}$ as $a(q^t) = \sum_{i=1}^{2g} \a_i^t.$ For $1 \le n \le \ell$, we have the following.
$$a(q^t) = \left\{
\begin{array}{ll}
0 & \text{if $m \nmid t$,}\\[3pt]
-\ff_{r,n} & \text{if $t=mr$ and $a^2-4b = \square$ in $\F_{q^m}$},\\[5pt]
(-1)^{r-1}\ff_{r,n} &\text{if $t=mr$ and $a^2-4b\ne \square$ in $\F_{q^m}$,}
\end{array} \right.$$
where all the notations be the same as in Theorem \ref{thmA}.

If $\ell \mid t$, then $a(q^t) = \left\{
\begin{array}{ll}
0 & \text{if $m \nmid t$,}\\[3pt]
-\ff_{r,\ell} & \text{if $t=mr$ and $a^2-4b = \square$ in $\F_{q^m}$},\\[5pt]
(-1)^{r-1}\ff_{r,\ell}& \text{if $t=mr$ and $a^2-4b \ne \square$ in $\F_{q^m}$:}
\end{array}
\right.$
that is, $n=\ell$.
\end{thm}

\begin{proof}
Let $K= k(\sqrt[\ell]{T^2+aT+b})$ be a function field which corresponds to \eqref{he} and $K_t = K\F_{q^t}$ be a constant field extension of $K$. Then the number of $\F_{q^t}$-points $\#HE_{a,b}(\F_{q^t})$ of $HE_{a,b}$ can be interpreted as the number of rational places of $K_t$: that is, the number of places of $K_t$ whose degrees are one.

We first consider the case where $m \mid t$. Let $t=m$ (i.e., $r=1$).
By Lemma \ref{hyperelliptic_ff}, we can easily find out that $\#HE_{a,b}(\F_{q^m}) = \left\{
\begin{array}{ll}
\ell\mathfrak{n}_{1,s}+3 & \text{if $a^2-4b=\square$ in $\F_{q^m}^\times$,}\\[5pt]
\ell\mathfrak{n}_{1,s}+1 & \text{if $a^2-4b\ne\square$ in $\F_{q^m}^\times$,}
\end{array}\right.$ where $\mathfrak{n}_{1,s}$ is the number of rational places of $\F_{q^m}(T)$ which splits completely in $K_m$. 
By Kummer's Theorem \cite[Theorem 3.3.7]{St}, for $c\in \F_{q^m}$ such that $a \ne c$, $b \ne c$, if
\begin{equation}\label{A}
y^\ell = T^2+aT+b \pmod {T-c}
\end{equation}
has roots over $\F_{q^m}$ then such $T-c$ is a place of $\F_{q^m}(T)$ which splits completely in $K_m$. We note that \eqref{A} can be interpreted as if there exists $c \in \F_{q^m}$ such that $c^2+ac+b = d^\ell$ for some $d \in \F_{q^m}^\times$. Thus, $\mathfrak{n}_{1,s}$ is equal to the number of the roots of $x^2+ax+b = d^\ell \text{ in } \F_{q^m};$ hence, $\mathfrak{n}_{1,s}$ becomes
$\mathfrak{n}_{1,s} = 2\mathfrak{m}_{1,\k,1}+\mathfrak{m}_{1,\k,2}$,
where $\kappa:= a^2-4b \in \F_{q}^\times$, and
\begin{align*}
\mathfrak{m}_{1,\kappa, 1} &:= \frac{|\{ d \in \F_{q^m}^\times \mid \kappa+4d^\ell = \eta^2 \text{ for some $\eta \in \F_{q^m}^\times$}\}| }{\ell},\\[3pt]
\mathfrak{m}_{1,\kappa,2} &:= \frac{|\{ d \in \F_{q^m}^\times \mid \kappa+4d^\ell = 0 \text{ in $\F_{q^m}$}\}|}{\ell}.
\end{align*}

We now compute $\mathfrak{m}_{1,\k,1}$ by using the Jacobi sum. For a fixed $\k \in \F_{q^m}^\times$, let $\mathcal{N}_{1,\kappa}$ be the number of solutions $(d, \eta) \in \F_{q^m}^2$ of $4d^\ell - \eta^2 = -\kappa$. By Lemma \ref{number}, we get
\allowdisplaybreaks
\begin{align*}
\mathcal{N}_{1,\kappa}
&= q^m + \sum_{j_1=1}^{\ell-1}\l_1^{j_1}((-\kappa)4^{-1})\l_2(\kappa)J(\l_1^{j_1}, \l_2),
\end{align*}
where $\l_1$ (resp. $\l_2$) is a multiplicative character of $\F_{q^m}$ of order $\ell$ (resp. 2) with $\l_1(0) = \l_2(0)=0$. For $c \in \F_{q^m}^\times$, we note that $\l_1(c) =
\zeta^n$ for $1 \le n \le \ell$ and $\l_2(c)=\pm 1$, where $\zeta:=e^{\frac{2\pi\imath}{\ell}}$. By the definition of Jacobi sum, we let
$J(\l_1,\l_2) := \sum_{i=1}^\ell a_{1,i}\z^i \quad \text{with }a_{1,i} \in \Z.$ Then we can easily find out that
$$J(\l_1^{j_1}, \l_2) = \sum_{c_1+c_2=1}\l_1^{j_1}(c_1)\l_2(c_2) = \sum_{i=1}^{\ell}a_{1,i}\zeta^{j_1i}.$$ 
Suppose that $\l_1((-\k)4^{-1}) = \z^n$ with $1 \le n \le \ell$. Then we obtain
\allowdisplaybreaks
\begin{align*}
\mathcal{N}_{1,\k} &= q^m + \sum_{j_1=1}^{\ell-1}\l_1^{j_1}((-\kappa)4^{-1})\l_2(\kappa)J(\l_1^{j_1},\l_2) \\[2pt]
&=\left\{
\begin{array}{ll}
q^m + \sum_{j_1=1}^{\ell-1}\l_1^{j_1}((-\kappa)4^{-1})J(\l_1^{j_1},\l_2) & \text{if $a^2-4b = \square$ in $\F_{q^m}^\times$,}\\[7pt]
q^m - \sum_{j_1=1}^{\ell-1}\l_1^{j_1}((-\kappa)4^{-1})J(\l_1^{j_1},\l_2) & \text{if $a^2-4b \ne \square$ in $\F_{q^m}^\times$,}
\end{array}\right.\\
&=q^m\pm\sum_{j_1=1}^{\ell-1}\z^{nj_1}\left(\sum_{i=1}^\ell a_{1,i}\z^{ij_1} \right) \\[2pt]
&= q^m\pm\sum_{j_1=1}^{\ell-1}\z^{nj_1}\left(a_{1,1}\z^{j_1}+a_{1,2}\z^{2j_1} + \cdots + a_{1,\ell-1}\z^{(\ell-1)j_1}+ a_{1,\ell}\z^{\ell j_1} \right)\\[3pt]
&=q^m\pm\sum_{i=1}^\ell a_{1,i}\left(\sum_{j_1=1}^{\ell-1}\z^{j_1(n+i)} \right) = q^m\pm((\ell-1)a_{1,\ell-n}-\sum_{\substack{i=1\\i\ne \ell-n}}^\ell a_i) = q^m \pm(\ell a_{1,\ell-n}-\sum_{i=1}^\ell a_{1,i}),
\end{align*}
where $a_{1,0}:=a_{1,\ell}$.

For simplicity, let $\ff_{1,n}:=\ell a_{1,\ell-n}-\sum_{i=1}^\ell a_{1,i}$. Then we get
\begin{align*}
\mathcal{N}_{1,\k} &= \# \{(d,\eta) \in \F_{q^m}^2 \mid 4d^\ell-\eta^2 = -\k \} =
\left\{
\begin{array}{ll}
q^m + \mathcal{F}_{1,n} & \text{if $a^2-4b = \square$ in $\F_{q^m}^\times$,}\\[5pt]
q^m - \mathcal{F}_{1,n} & \text{if $a^2-4b \ne \square$ in $\F_{q^m}^\times$,}\end{array}
\right.
\end{align*}
where $\l_1(-4^{-1}(a^2-4b))=\z^n$. Suppose $a^2-4b$ is a square in $\F_{q^m}^\times$: that is, $\kappa$ is a square. Since we only need the number of $d \in \F_{q^m}^\times$, we need to subtract 2 from $\mathcal{N}_{1,\k}$; this is because $(0,\pm \eta)$ are the roots of $4d^\ell-\eta^2=-\k$. Furthermore, for the same $d$, there are two possibilities of $\eta$: namely, $\pm \eta$. Thus we get
$$\mathfrak{m}_{1,\k,1} =\left\{
\begin{array}{ll}
\mathlarger{\frac{1}{\ell}\cdot \frac{q^m+\mathcal{F}_{1,n}-2}{2}} & \text{if $\l_1(-4^{-1}\k) = \zeta^n$ with $1 \le n < \ell$},\\[9pt]
\mathlarger{\frac{1}{\ell}\cdot \frac{q^m+\mathcal{F}_{1,n}-2-\ell}{2}} & \text{if $\l_1(-4^{-1}\k) = 1$}.
\end{array}\right.$$

For the case where $a^2-4b$ is a non-square in $\F_{q^m}^\times$, we have the following:
$$\mathfrak{m}_{1,\k,1} =\left\{
\begin{array}{ll}
\mathlarger{\frac{1}{\ell}\cdot \frac{q^m-\mathcal{F}_{1,n}}{2}} & \text{if $\l_1(-4^{-1}\k) = \zeta^n$ with $n \ne  \ell$},\\[9pt]
\mathlarger{\frac{1}{\ell}\cdot \frac{q^m-\mathcal{F}_{1,n}-\ell}{2}} & \text{if $\l_1(-4^{-1}\k) = \z^\ell = 1$}.
\end{array}\right.$$

For estimation of $\mathfrak{m}_{1,\k,2}$, we divide into 2 cases: namely, $\l_1(-4^{-1}\k) \ne 1$ \mbox{and $\l_1(-4^{-1}\k) = 1$.} If $\l_1(-4^{-1}\k) \ne 1$, then there is no element in $\F_{q^m}^\times$ which makes $\ell$th root. On the otherhand, if $\l_1(-4^{-1}\k) = 1$, then there exists $\epsilon \in \F_{q^m}^\times$ such that $(-\k^2)4^{-1} = \epsilon^\ell$. Therefore, we get
$$\mathfrak{m}_{1,\k,2} =\left\{
\begin{array}{ll}
\mathlarger{\frac{1}{\ell}\cdot 0 = 0} & \text{if $\l_1(-4^{-1}\k) \ne 1$}, \\[9pt]
\mathlarger{\frac{1}{\ell}\cdot \ell = 1} & \text{if $\l_1(-4^{-1}\k) = 1$}.
\end{array}\right.$$
Thus, if $a^2-4b$ is a square in $\F_{q^m}^\times$, we have the desired result as follows:
\begin{align*}
&\ell\mathfrak{n}_{1,s}+3 = \ell(2\mathfrak{m}_{1,\k,1}+\mathfrak{m}_{1,\k,2})+3\nonumber\\[5pt]
&= \left\{
\begin{array}{ll}
\mathlarger{\ell \cdot  \left(2\cdot \frac{1}{\ell}\cdot \frac{q+\mathcal{F}_{1,n}-2}{2} + 0\right)+3} & \text{if $\l_1(-4^{-1}\k) = \zeta^n$ with $1 \le n < \ell$},\\[9pt]
\mathlarger{\ell\cdot  \left(2\cdot \frac{1}{\ell}\cdot \frac{q+\mathcal{F}_{1,n}-2-\ell}{2} + 1 \right)+3} & \text{if $\l_1(-4^{-1}\k) = \z^n$ with $n=\ell$,}
\end{array}\right.\\[5pt]
&= q^m + 1 + \ff_{1,n}.
\end{align*}
If $a^2-4b$ is a non-square in $\F_{q^m}^\times$,
\allowdisplaybreaks
\begin{align*}
&\ell\mathfrak{n}_{1,s}+1 = \ell(2\mathfrak{m}_{1,\k,1}+\mathfrak{m}_{1,\k,2})+1\nonumber\\[5pt]
&= \left\{
\begin{array}{ll}
\mathlarger{\ell \cdot  \left(2\cdot \frac{1}{\ell}\cdot \frac{q-\mathcal{F}_{1,n}}{2}\right)+1} = q-\mathcal{F}_{1,n}+1 & \text{if $\l_1(-4^{-1}\k) = \zeta^n$ with $n \ne \ell$},\\[11pt]
\mathlarger{\ell\cdot  \left(2\cdot \frac{1}{\ell}\cdot \frac{q-\mathcal{F}_{1,n}-\ell}{2}+1\right)+1} = q-\mathcal{F}_{1,n}+1 & \text{if $\l_1(-4^{-1}\k) = 1$}.
\end{array}\right.
\end{align*}
Using the relation between $a(q^m)$ and $\#HE_{a,b}(\F_{q^m})$, we get the desired result.

Now, let $t=mr$ with $r \ge 2$. By \cite[Corollary 5.16]{St}, we have $a(q^t) = q^t+1-\#HE_{a,b}(\F_{q^t})$.
Let $\mathcal{N}_{r,\kappa}$ be the number of solutions $(d, \eta) \in \F_{q^t}^2$ of $4d^\ell - \eta^2 = -\k$ with $\k=a^2-4b$. Let $\l_{r,1}$ (resp. $\l_{r,2}$) be a multiplicative character of $\F_{q^t}$ of order $\ell$ (resp. 2) such that $\l_{r,i} = \l_i \circ N_{\F_{q^{mr}}/\F_{q^m}}$ and $\l_{r,i}(0) =0$ ($i=1,2$), where $N_{\F_{q^{mr}}/\F_{q^m}}$ denotes the norm of $\F_{q^{mr}}$ over $\F_{q^m}$. Then by Lemma \ref{jacobi_lift}, we obtain
\begin{equation*}
J(\l_{r,1},\l_{r,2}) = (-1)^{(r-1)(2-1)}J(\l_1,\l_2)^r.
\end{equation*}
Letting $(-1)^{(r-1)}J(\l_1,\l_2)^r := \sum_{i=1}^\ell a_{r,i}\z^i$, we get
\allowdisplaybreaks
\begin{align*}
\mathcal{N}_{r,\k} &= q^t+\sum_{j_1=1}^{\ell-1}\l_{r,1}^{j_1}(-4^{-1}\k)\l_{r,2}(\k)J(\l_{r,1}^{j_1},\l_{r,2})\nonumber\\[3pt]
&= \left\{
\begin{array}{ll}
q^t+ \sum_{j_1=1}^{\ell-1}\l_{r,1}^{j_1}(-4^{-1}\k)J(\l_{r,1}^{j_1},\l_{r,2}) & \text{if $a^2-4b=\square$ in $\F_{q^m}^\times$},\\[7pt]
q^t+(-1)^{r}\sum_{j_1=1}^{\ell-1}\l_{r,1}^{j_1}(-4^{-1}\k)J(\l_{r,1}^{j_1},\l_{r,2}) & \text{if $a^2-4b\ne \square$ in $\F_{q^m}^\times$},
\end{array}
\right.\\
&=\left\{
\begin{array}{ll}
q^t+\ell a_{r,\overline{\ell-rn}} - \sum_{i=1}^\ell a_{r,i} & \text{if $a^2-4b = \square$ in $\F_{q^m}^\times$,} \\[7pt]
q^t+(-1)^{r}\left(\ell a_{r,\overline{\ell-rn}} - \sum_{i=1}^\ell a_{r,i}\right) & \text{if $a^2-4b \ne \square$ in $\F_{q^m}^\times$,}
\end{array}
\right.
\end{align*}
where $\l_{r,1}(-4^{-1}\k) = \z^n$. The second equality follows from the fact $\l_{r,2}(\k) = \l_2(N_{\F_{q^{mr}}/\F_{q^m}}(\k)) = \l_2(\k^r)=(\l_2(\k))^r$. Also, we note that for $c\in \F_{q^m}^\times$, if $\l_1(c) = \z^n$, then $\l_{r,1}(c) = \z^{rn}$; therefore, we get the desired result.

We now assume that $m \nmid t$. By Lemma \ref{hyperelliptic_ff}, we note that the infinite place $\mathfrak{p}_\infty$ of $\F_{q^t}(T)$ is totally ramified in $K_t$. Let $T-c$ be a rational place of $\F_{q^t}(T)$ with $a \ne c$, $b\ne c$. Then the ramification behavior of $T-c$ depends on the factorization of $y^\ell \equiv T^2+aT+b \pmod{T-c}$; this is equivalent to
$y^\ell-(c^2+ac+b)$ over $\F_{q^t}$. Since $q^t \not \equiv 1 \pmod \ell$, the equation $y^\ell-(c^2+ac+b)$ always has the unique root in $\F_{q^t}$; this implies that the rational place $T-c$ splits partially in $K_t$. Thus, the number of rational places of $K_t$ is
$$\left\{
\begin{array}{ll}
\ell\times 0 + 3 + q^t-2= q^t+1 & \text{if $a^2-4b=\square$ in $\F_{q^m}^\times$}\\[3pt]
\ell\times 0+1+q^t=q^t+1 & \text{if $a^2-4b \ne \square$ in $\F_{q^m}^\times$;}
\end{array}\right.$$
thus, we get the desired result.

Finally, assume that $\ell \mid r$. If $m \nmid t$, then the fact $a(q^t)=0$ is always zero follows immediately. We now assume that $m \mid t$. As before, letting $\l_1(-4^{-1}\k)=\z^n$, then
$$\ff_{r,n} = \ell a_{r,\overline{\ell-rn}}-\sum_{i=1}^\ell a_{r,i} = \ell a_{r,0} - \sum_{i=1}^\ell a_{r,i} = \ff_{r,\ell};$$
this implies that the value of $a(q^t)$ does not depend on $\l_1(-4^{-1}\k)=\z^n$.
\end{proof}

Let all the notations be the same as in Theorem \ref{thm1}. In the following remark, we note that the value of $a(q^t)$ does not depend on the choice of multiplicative character $\l_1$ of $\F_{q^m}$ of order $\ell$ even though there are $(\ell-1)$ distant $\l_1$. We give the explanation as follows.

\begin{rem}
Let $g$ and $g'$ be distinct generators of $\F_q^\times$ and $\l_1(g) = \l_1'(g') = \z$, where $\l_1$ and $\l_1'$ are distinct two multiplicative characters of $\F_{q^m}$ of order $\ell$ and $\z = e^{\frac{2\pi \imath}{\ell}}$. Then there exists some $\delta$ ($2 \le \delta \le \ell-1)$ such that $\l_1'(g)=\z^\d$. Let $J(\l_1,\l_2) = \sum_{i=1}^\ell a_{1,i}\z^i$ with a quadratic multiplicative character $\l_2$ of $\F_{q^m}$. Then we obtain
\begin{align*}
J(\l_1',\l_2) &=\sum_{c_1+c_2=1}\l_1'(c_1)\l_2(c_2) = \sum_{c_1+c_2=1}\l_1^\d(c_1)\l_2(c_2) = \sum_{i=1}^\ell a_{1,i}\z^{i\d}.
\end{align*}
Therefore, we get
\begin{align*}
(-1)^{r-1}J(\l_1',\l_2)^r=(-1)^{r-1}\left(\sum_{i=1}^\ell a_{1,i}\z^{i\d}\right)^r= \sum_{i=1}^\ell a_{r,\overline{i\d^{-1}}}\z^i,
\end{align*}
where we use the notation $(-1)^{r-1}J(\l_1,\l_2)^r = \sum_{i=1}^\ell a_{r,i}\z^i$ and $\d^{-1}$ denotes the inverse element of $\d$ in $\Z/\ell\Z$. We let
$$(-1)^{r-1}J(\l_1',\l_2)^r = \sum_{i=1}^\ell a'_{r,i}\z^i:$$
that is, $a'_{r,i} = a_{r,\overline{i\d^{-1}}}$. Now, assume that $\l_1(-4^{-1}(a-b)^2)=\z^n$; then, $\l_1'(-4^{-1}(a-b)^2)=\z^{n\d}$. Using the relation $(\ell-rn\d)\d^{-1} \equiv -rn \equiv \ell - rn \pmod \ell$, we have
\begin{align*}
\ff'_{r,n\d} = \ell a'_{r,\overline{\ell-r(n\d)}} - \sum_{i=1}^\ell a'_{r,i} = \ell a_{r,\overline{\ell-rn}}-\sum_{i=1}^\ell a_{r,i} = \ff_{r,n};
\end{align*}
hence, $a(q^t)$ does not depend on the choice of $\l_1$.
\end{rem}

\vskip 0.3cm

We provide the following lemma, which gives the closed formula of $S_{rm}$, where $m$ is even and $r \in \Z^+$. This result plays an important role in the proof of Theorem \ref{thmB}.

\begin{lem}\label{S2}
Let $m$ be the order of odd prime power $q$ in $\F_\ell^\times$. If $m$ is even, then
$$S_{rm} = (-1)^{r-1}(\ell-1)q^{\frac{rm}{2}},$$
where $r \in \Z^+$ and $S_{rm}$ is defined in Lemma \ref{c}.
\end{lem}

\begin{proof}
From the assumption that $m$ is even, we have $q^{\frac{m}{2}} \equiv -1 \pmod \ell$ and $q^m \equiv 1 \pmod \ell$. Let $J(\l_1, \l_2)$ be a Jacobi sum, where $\l_1$ (resp. $\l_2$) is a multiplicative character of $\F_{q^m}$ of order $\ell$ (resp. 2). Then, we have
\begin{align}\label{real}
J(\lambda_1, \lambda_2) =J(\lambda_1^{q^{\frac{m}{2}}}, \lambda_2^{q^{\frac{m}{2}}})= J(\lambda_1^{-1}, \lambda_2)=\overline{J(\lambda_1, \lambda_2)};
\end{align}
the first identity holds by the Frobenius automorphism for $\F_{q^m}/\F_q$ and the second one is true because $q^{\frac{m}{2}} \equiv -1 \pmod \ell$. Therefore, \eqref{real} shows that $J(\lambda_1, \lambda_2)$ is real and its value is equal to $q^{\frac{m}{2}}$ or $-q^{\frac{m}{2}}$. Since $\lambda_1$, $\lambda_2$, and $\lambda_1\lambda_2$ are not trivial, we have the following identity.
$$J(\lambda_1, \lambda_2)=g(\lambda_1)g(\lambda_2)/g(\lambda_1\lambda_2),$$
where $g(\l)$ is the Gauss sum.
On the other hand,
\begin{align*}
 g(\chi)^\ell=\bigg(\sum_t \chi(t)\zeta_p^{\text{Tr}(t)}\bigg)^\ell \equiv \sum_t \chi(t)^\ell\zeta_p^{\text{Tr}(\ell t)} \equiv  \sum_{t\neq 0} \zeta_p^{\text{Tr}(\ell t)}\equiv -1 \pmod \l
\end{align*}
in $\Z[\zeta,\zeta_p]$; this implies that $J(\lambda_1, \lambda_2)^\ell\equiv -1 \pmod \l$. If $J(\lambda_1, \lambda_2)=-q^{\frac{m}{2}}$, then $J(\lambda_1, \lambda_2)^\ell = -q^\ell \equiv -(-1)^\ell \equiv 1 \pmod \l$, which is a contradiction. Hence $J(\lambda_1, \lambda_2)=q^{\frac{m}{2}}$.

Let $J(\l_1,\l_2)=\sum_{i=1}^\ell a_{1,i}\z^i$, where $\z$ is the $\ell$th root of unity. Since $J(\l_1, \l_2)$ is real, we get $a_{1,1} = a_{1,2} = \cdots = a_{1,\ell-1}$; hence $J(\l_1,\l_2) = a_{1,\ell}-a_{1,1}=q^{\frac{m}{2}}$. By the fact that $a^2-4b$ is always a square in $\F_{q^m}^\times$ and $(-1)^{r-1}J(\l_1,\l_2)^r = \sum_{i=1}^\ell a_{r,i}\z^i =  (-1)^{r-1}q^{\frac{rm}{2}}$, we obtain the following:
$$S_{rm} = \ell a_{r,\ell} - \sum_{i=1}^\ell a_{r,i} =\ell(-1)^{r-1}q^{\frac{rm}{2}} - ((-1)^{r-1}q^{\frac{rm}{2}}) = (-1)^{r-1}(\ell-1)q^{\frac{rm}{2}}.$$
\end{proof}

We now give the proof of Theorem \ref{thmB}.

\textbf{Proof of Theorem \ref{thmB}}

Let $m \le \ell-1$ and even. By Theorem \ref{thmA}, it suffices to find the value of $c_m, c_{2m}, \cdots, c_{rm}$, where $r:= \lfloor \frac{g}{m} \rfloor$. We note that for $i \nmid m$, the value $c_i=0$ since $S_i=0$. We claim the following:
\begin{equation}\label{cm}
c_{rm} = \binom{\frac{\ell-1}{m}}{r}q^{\frac{rm}{2}}.
\end{equation}
Assume that $m \mid g$: that is, $\frac{\ell-1}{2} = rm$. Then the $L$-polynomial of $K$ is

\begin{align*}
L_K(u) &= \sum_{i=0}^{r-1}(q^{(r-i)m}u^{(2r-i)m}+u^{im})c_{im} + c_{rm}u^{rm}\\
&=\sum_{i=0}^{r-1}\binom{\frac{\ell-1}{m}}{i}(q^{\frac{\ell-1-im}{2}}u^{(2r-i)m}+q^{\frac{im}{2}}u^{im})+\binom{\frac{\ell-1}{m}}{r}q^{\frac{rm}{2}};
\end{align*}
we use \eqref{cm} for the second equality. Using the fact that $L_K(1)$ is the class number $h_K$ of $K$, we further get
\begin{align*}
h_K = L_K(1) &= \sum_{i=0}^{r-1}\binom{2r}{i}(q^{\frac{(2r-i)m}{2}}+q^{\frac{im}{2}}) + \binom{2r}{r}q^{\frac{rm}{2}}=\sum_{i=0}^{2r}\binom{2r}{i}q^{\frac{im}{2}} = (q^{\frac{m}{2}}+1)^{2r};
\end{align*}
we use $\frac{\ell-im-1}{2}=\frac{(2r-i)m}{2}$ and $\binom{2r}{i} = \binom{\frac{\ell-1}{m}}{2r-i}$ for the third equality. We now suppose $m \nmid g$. In this case, we have
$$r = \bigg\lfloor \frac{g}{m} \bigg\rfloor = \bigg\lfloor \frac{\frac{\ell-1}{2}}{m} \bigg\rfloor = \bigg\lfloor \frac{\frac{\ell-1}{m}}{2} \bigg\rfloor \quad \Leftrightarrow \quad 2r+1 = \frac{\ell-1}{m}.$$
Using $2r+1 = \frac{\ell-1}{m}$ and a similar reasoning as in $m \mid g$, we have the desired result.

We now prove \eqref{cm}. For the proof of \eqref{cm}, we need the following:
\begin{equation}\label{claim}
(-1)^{s-1}(s-1)!m^{s-1} + \sum_{j=1}^{s-1}\bigg((-1)^{s-j-1}\frac{(s-1)!}{j!}m^{s-j-1}\prod_{i=0}^{j-1}(\ell-im-1) \bigg) = \prod_{i=1}^{s-1}(\ell-im-1),
\end{equation}
where $s \in \Z_{\ge 2}$. We use Mathematical Induction on $s$. If $s=2$, then both left-hand side and right-hand side of \eqref{claim} are $\ell-m+1$; thus, \eqref{claim} holds. We now suppose \eqref{claim} is true for $s$. Then we get
\begin{align*}
&(-1)^{s}(s)!m^{s} + \sum_{j=1}^{s}\bigg((-1)^{s-j}\frac{s!}{j!}m^{s-j}\prod_{i=0}^{j-1}(\ell-im-1) \bigg)\\
&=(-sm)\bigg((-1)^{s-1}(s-1)!m^{s-1} + \sum_{j=1}^{s-1}\bigg((-1)^{s-j-1}\frac{(s-1)!}{j!}m^{s-j-1}\prod_{i=0}^{j-1}(\ell-im-1) \bigg)\bigg) + \prod_{i=0}^{s-1}(\ell-im-1)\\
&=(-sm)\prod_{i=1}^{s-1}(\ell-im-1)+\prod_{i=0}^{s-1}(\ell-im-1) = \prod_{i=1}^s(\ell-im-1);
\end{align*}
hence, \eqref{claim} holds. We now show \eqref{cm} by Mathematical Induction on $r$. By Lemma \ref{S2}, $S_{rm} = (-1)^{r-1}(\ell-1)q^{\frac{rm}{2}}$. For the case $r=1$ holds immediately. Under the assumption that \eqref{cm} holds for $r$, we have the following:
\allowdisplaybreaks
\begin{align*}
c_{(r+1)m}&=\frac{S_{(r+1)m}+S_{rm}c_m+\cdots + S_mc_{rm}}{(r+1)m}\\[5pt]
&=\frac{(-1)^r(\ell-1)q^{\frac{(r+1)m}{2}}+ (-1)^{r-1}(\ell-1)q^{\frac{rm}{2}}\times \frac{\ell-1}{m}q^{\frac{m}{2}}+ \cdots + (\ell-1)q^{\frac{m}{2}}\times \frac{\prod_{i=0}^{r-1}(\ell-im-1)}{r!m^r}q^{\frac{rm}{2}}}{(r+1)m}\\[5pt]
&=\frac{(\ell-1)q^{\frac{(r+1)m}{2}}}{(r+1)!m^{r+1}}\bigg((-1)^rr!m^r + (-1)^{r-1}r!m^{r-1} + (-1)^{r-2}\frac{r!}{2!}m^{r-2}(\ell-1)(\ell-m-1)+ \\
& \hskip 3.3cm \cdots + \prod_{i=0}^{r-1}(\ell-im-1) \bigg)\\[5pt]
&=\frac{(\ell-1)q^{\frac{(r+1)m}{2}}}\prod_{i=1}^r(\ell-im-1) = \frac{\prod_{i=0}^r(\ell-im-1)}{(r+1)!m^{r+1}}q^{\frac{(r+1)m}{2}} = \binom{\frac{\ell-1}{m}}{r+1}q^{\frac
{(r+1)m}{2}};
\end{align*}
thus, \eqref{cm} holds.

We now let $m = \frac{\ell-1}{2}$ and odd: that is, $m$ is equal to $g$. By Theorem \ref{thmA} and Theorem \ref{thm1}, we obtain the following:
$$S_g = S_m = \left\{
\begin{array}{ll}
-\ff_{1,n} & \text{if $a^2-4b = \square$ in $\F_{q^m}^\times$,}\\[3pt]
\ff_{1,n} & \text{if $a^2-4b \ne \square$ in $\F_{q^m}^\times$,}
\end{array}
\right.\text{ and $S_i=0$ for $1 \le i \le m-1$};$$
hence,
\begin{align*}
L_K(u) &= c_0(u^0+q^{g-0}u^{2g})+c_gu^g = q^gu^{2g} + 1 + \frac{S_m}{m}u^g\\[4pt]
&= \left\{
\begin{array}{ll}
q^gu^{2g}+ 1 -\frac{\ff_{1,n}}{m}u^g & \text{if $a^2-4b = \square$ in $\F_{q^m}^\times$,}\\[3pt]
q^gu^{2g} + 1 +\frac{\ff_{1,n}}{m}u^g & \text{if $a^2-4b \ne \square$ in $\F_{q^m}^\times$,}
\end{array}
\right.
\end{align*}
and $h_K = q^m + 1 + \frac{S_m}{m}$. Therefore, it is sufficient to show that $\ff_{1,n} = \ell a_{1,\ell}+1$. We have
\begin{equation}\label{lem_7}
\sum_{i=1}^\ell a_{1,i} = -1;
\end{equation}
this holds because of the following reasoning. As $\l_2$ has order 2, the value of $\l_2(c_2)$ is $\pm 1$ with $J(\l_1,\l_2) = \sum_{\substack{c_1+c_2=1\\c_1, c_2 \in \F_{q^m}^\times}} \l_1(c_1)\l_2(c_2)$. Therefore, $\sum_{i=1}^\ell a_{1,i}$ is determined by subtracting the number of quadratic non-residue from that of quadratic residue. In general, the numbers of quadratic residue and quadratic non-residue are the same; however, since Jacobi sum does not allow the value $c_1=0$, \eqref{lem_7} is true. Also, we note that for every element $c \in \F_q^\times$, the value of $\l_1(c)$ is always 1: that is, $\z^\ell$. Therefore,
$\ff_{1,n} = \ell a_{1,\overline{\ell-\ell}} - \sum_{i=1}^\ell a_{1,i} = \ell a_{1,\ell}+1$.
\qed


\section{Applications for the $L$-polynomial of $k(\sqrt[\ell]{T^2+aT+b})$\\ : the average class numbers of $\HE_5(\F_q)$ and $\HE_7(\F_q)$}

In this section, we give explicit formulas for the average class numbers of $\HE_5(\F_q)$ and $\HE_7(\F_q)$ in Theorem \ref{thmC} and Theorem \ref{thmD}, respectively: that is,
\begin{align*}
&\HE_5(\F_{q}):=\{\F_q(T,\sqrt[5]{T^2+aT+b}) \mid a, b \in \F_q,~ a^2-4b \ne 0\} \quad \text{and}\\[5pt]
&\HE_7(\F_{q}):=\{\F_q(T,\sqrt[7]{T^2+aT+b}) \mid a, b \in \F_q,~ a^2-4b \ne 0\}.
\end{align*}

For the proofs, we introduce two lemmas, which play a crucial role for the proofs of Theorems \ref{thmC} and \ref{thmD}. Lemma \ref{cor1} states the average value on the trace of Frobenius of (hyper)elliptic curves $HE_{a,b}$ over $\F_{q^t}$ with $m=1$. Lemma \ref{lem_5} gives the relation between the coefficients of Jacobi sum $J(\l_1,\l_2) = \sum_{i=1}^\ell a_i\z^i$ with $\ell=5,7$ and $\z=e^{\frac{2\pi \imath}{\ell}}$, where $\l_1$ (resp. $\l_2$) is a multiplicative character of $\F_{q}$ of order $\ell$ (resp. order 2).

\vskip 0.3cm

\begin{lem}\label{cor1}
Let $HE_{a,b}$ be a $($hyper$)$elliptic curve which is defined in \eqref{he} with $m=1$. Let all the notations be the same as in Theorem \ref{thmA}. Then the average value on $a(q^t)$ is
$$a(q^t)^{\av}=\left\{
\begin{array}{ll}
0& \text{if $\ell \nmid t$},\\[3pt]
-\ff_{t,\ell} & \text{if $\ell \mid t$ and $a^2-4b = \square$ in $\F_q^\times$},\\[3pt]
(-1)^{t-1}\ff_{t,\ell} & \text{if $\ell \mid t$ and $a^2-4b\ne \square$ in $\F_q^\times$.}
\end{array}\right.$$
\end{lem}

\begin{proof}
We divide into two cases: namely, whether $a^2-4b$ is a square or a non-square in $\F_q^\times$.

\textbf{Case 1.} $a^2-4b$ is square in $\F_q^\times$.

Let $t$ be a positive integer which is not divisible by $\ell$. Then all possible values of $a(q^t)$ is $-\ff_{t,n}$ depending on $\l_1(-4^{-1}(a^2-4b))=\z^n$. We note that among the $\frac{q(q-1)}{2}$ possible hyperelliptic curves $y^\ell=x^2+ax+b$, there are $\frac{q(q-1)}{2\ell}$ hyperelliptic curves such that $\l_1(-4^{-1}(a-b)^2) =\z^n$ for each $1 \le n \le \ell$. Therefore,
\begin{align*}
a(q^t)^{\av} &=\frac{2}{q(q-1)}\sum_{\substack{a,b \in \F_q\\a^2-4b=\square}}a(q^t) = \frac{2}{q(q-1)}\cdot \frac{q(q-1)}{2\ell} \cdot \sum_{i=1}^\ell(-\ff_{t,i})=0.
\end{align*}
For the case where $\ell \mid t$, the value of $a(q^t)$ is always $-\ff_{t,\ell}$; hence, the result follows immediately.

\textbf{Case 2.} $a^2-4b$ is a non-square in $\F_q^\times$.

We claim that there are $\frac{q(q-1)}{2\ell}$ hyperelliptic curves such that $\l_1(-4^{-1}(a^2-4b))=\z^n$ for each $1 \le n \le \ell$, where $a^2-4b$ is a non-square over $\F_q^\times$ and $\l_1$ is a multiplicative character of $\F_q$ of order $\ell$. Since $\l_1(-4^{-1}(a^2-4b))=\l_1(-4^{-1}) \l_1(a^2-4b)=\z^\epsilon \cdot \z^m$ with $\l_1(-4^{-1})=\z^\epsilon$, it is enough to show that $m$ takes all the values among 1 through $\ell$ with the same multiplicity $\frac{(q-1)}{2\ell}$. We note that there are $\frac{q-1}{2}$ numbers of $c:=a^2-4b \in \F_q^\times$ and the value $m$ appears repetitively with period $\ell$. Furthermore, the value $m$ should be odd since $c$ is non-square. Therefore, it is sufficient to consider the first $\ell$ numbers of odd integers, namely: $\{1, 3, 5, \cdots, 2\ell-1\}$.
\allowdisplaybreaks
\begin{align*}
&\{1, 3, 5, \cdots, 2\ell-1\} \\
&= \left\{1, 2\times 1+1, 2\times 2+1, \cdots, 2\times \frac{\ell-1}{2}+1, 2\times \frac{\ell+1}{2}+1, \cdots,  2\times(\ell-1)+1\right\}\\[3pt]
&\equiv \{1,3,5, \cdots, \ell-2, 0, 2, 4, \cdots, \ell-1\} \pmod \ell;
\end{align*}
we use the fact that $\ell$ is an odd prime. Since the first appearing $\ell$ numbers of $m$ are distinct, we prove the claim.

Using a similar reasoning as in the proof of Case 1, we obtain the following. We use the fact that there are $\frac{q(q-1)}{2}$ numbers of monic quadratic irreducible \hbox{polynomial in $\F_q[x]$.}
\allowdisplaybreaks
\begin{align*}
a(q^t)^{\av} = \frac{2}{q(q-1)}\sum_{\substack{a,b\in \F_q\\
a^2-4b \ne \square}} a(q^t) &=\left\{
\begin{array}{ll}
\frac{2}{q(q-1)}\cdot \frac{q(q-1)}{2\ell} \cdot \sum_{i=1}^\ell((-1)^{t+1}\ff_{t,i})=0 & \text{if $\ell \nmid t$}\\[10pt]
\frac{2}{q(q-1)}\cdot \frac{q(q-1)}{2} \cdot ((-1)^{t+1}\ff_{t,\ell})= (-1)^{t-1}\ff_{t,\ell} & \text{if $\ell \mid t$.}
\end{array} \right.
\end{align*}
\end{proof}


\begin{lem}\label{lem_5}
Let $\z$ be the $\ell$th root of unity and $J(\l_1, \l_2) = \sum_{i=1}^\ell a_i\z^i$ with $\ell = 5, 7$, where $a_i \in \Z$, $\l_1$ $($resp. $\l_2)$ is a multiplicative character of $\F_q$ of order $\ell$ $($resp. $2)$. Then we have the following:
\begin{equation}\label{q}
q=\sum_{i=1}^\ell a_i^2 - \frac{2}{\ell-1}\sum_{
\substack{1 \le i < j \le \ell}}
a_ia_j.
\end{equation}
For $\ell=7$, we further have the following.
$$\sum_{i=1}^7 (a_i^2) = \frac{6q+1}{7}, \quad \sum_{1\le i< j \le 7}^7a_ia_j = \frac{-3q+3}{7}, \quad \sum_{i=1}^7(a_i^3) - 3\sum_{1\le i <j <k\le 7}a_ia_ja_k = \frac{-9q+2}{7}.$$
\end{lem}

\begin{proof}
We first consider the case where $\ell=5$. Then $\z = \cos \frac{2\pi}{5}+\imath \sin \frac{2\pi}{5}$. Using the fact $|J(\l_1, \l_2)| = |\sum_{i=1}^5 a_i\z^i| =
\sqrt{q}$ \cite[p. 209, Theorem 5.22]{LN}, we obtain the following:
\allowdisplaybreaks
\begin{align*}
q= |J(\l_1, \l_2)|^2 = \bigg|\sum_{i=1}^5 a_i\z^i\bigg|^2 &= \sum_{i=1}^5a_i^2 + (2a_1a_2+2a_1a_5+2a_2a_3+2a_3a_4+2a_4a_5)\cos\frac{2\pi}{5}\\
&\hskip 1.5cm + (2a_1a_3+2a_1a_4+2a_2a_4+2a_2a_5+2a_3a_5)\cos\frac{4\pi}{5}\\[3pt]
&=\sum_{i=1}^5a_i^2 + (2a_1a_2+2a_1a_5+2a_2a_3+2a_3a_4+2a_4a_5)\frac{\sqrt{5}-1}{4}\\
&\hskip 1.5cm + (2a_1a_3+2a_1a_4+2a_2a_4+2a_2a_5+2a_3a_5)\frac{-\sqrt{5}-1}{4};
\end{align*}
thus, the result follows immediately.

We now assume that $\ell=7$. Then we have
\allowdisplaybreaks
\begin{align}\label{7}
q= |J(\l_1, \l_2)|^2 &= \sum_{i=1}^7a_i^2 + 2\cos\frac{2\pi}{7}\underbrace{(a_1a_2+a_2a_3+a_3a_4+a_4a_5+a_5a_6+a_6a_7+a_1a_7)}_{\mathcal{A}}\nonumber\\
&\hskip 1.5cm + 2\cos\frac{4\pi}{7}\underbrace{(a_1a_3+a_2a_4+a_3a_5+a_4a_6+a_1a_6+a_2a_7+a_5a_7)}_{\mathcal{B}}\nonumber\\[2pt]
&\hskip 1.5cm + 2\cos\frac{6\pi}{7}\underbrace{(a_1a_4+a_2a_5+a_3a_6+a_4a_7+a_1a_5+a_2a_6+a_3a_7)}_{\mathcal{C}}.
\end{align}
Let $r_1, r_2, r_3 \in \R$ be the roots of $r^3+r^2-2r-1=0$ with $r_1 >0$ and $r_2=r_1^2-2, r_3 < 0$. Then we note that $\z = e^{\frac{2\pi \imath}{7}} = \frac{r_1}{2} + \imath \sqrt{1 - \frac{r_1^2}{4}}$; thus, $2\sum_{j=1}^3\cos\frac{2j\pi}{7}=2\left(\frac{r_1}{2}+\frac{r_2}{2}+\frac{r_3}{2}\right) = -1$. Hence \eqref{7} becomes
\begin{align*}
q=\sum_{i=1}^7a_i^2 + r_1(\mathcal{A}-\mathcal{C})+r_2(\mathcal{B}-\mathcal{C})-\mathcal{C} & \quad \Leftrightarrow \quad r_1(\mathcal{A}-\mathcal{C})+r_2(\mathcal{B}-\mathcal{C}) \in \Z \\
& \quad \Leftrightarrow \quad r_1^2(\mathcal{B}-\mathcal{C}) + r_1(\mathcal{A}-\mathcal{C})-2(\mathcal{B}-\mathcal{C}) \in \Z,
\end{align*}
this only happens when $\mathcal{A}-\mathcal{C} = \mathcal{B}-\mathcal{C}=0$. Therefore,
$$q = \sum_{i=1}^7a_i^2 + 2\mathcal{A}\left(\cos\frac{2\pi}{7}+\cos\frac{4\pi}{7}+\cos\frac{6\pi}{7} \right) = \sum_{i=1}^7a_i^2-\mathcal{A} = \sum_{i=1}^7a_i^2 - \frac{1}{3}\sum_{1 \le i<j\le 7}a_ia_j.$$
The rest of equalities follow immediately by combining \eqref{lem_7} and \eqref{q}.
\end{proof}

We are now ready to prove Theorems \ref{thmC} and \ref{thmD}.

\begin{thm}\label{thmC}
Let $q$ be an odd prime power and $m$ be the order of $q$ in $\F_5^\times$. The average class number of $\HE_5(\F_q) = \{k(\sqrt[5]{T^2+aT+b}) \mid a,b \in \F_q, ~a^2-4b \ne 0\}$ is
$$(h_K)^{\av} = \left\{
\begin{array}{ll}
(q+1)^2& \text{if $m = 1, 2~~$ (i.e., $q \equiv \pm 1 \pmod 5$),}\\[5pt]
q^2+1  &\text{if $m=4~~$ (i.e., $q \not\equiv \pm 1 \pmod 5$).}
\end{array}
\right.$$

In fact, the average class number does not depend on whether $a^2-4b$ is a square or non-square in $\F_{q^m}^\times$.
\end{thm}

\begin{proof}
Let $K \in \HE_5(\F_{q})$ be a hyperelliptic function field. Since $\ell=5$, the possible $m$ is 1,2, or 4. By Theorem \ref{thmB}, it is sufficient to compute $m=1$.

Let the order of $q$ in $\F_5^\times$ be 1 and $L_K(u) = \sum_{i=0}^4c_iu^i \in \Z[u]$ the $L$-polynomial of $K$. By Lemma \ref{c}, the class number of $K$ can be computed as follows:
\begin{equation*}
h_K=(q^2+1)+(q+1)c_1+c_2 = (q^2+1) + (q+1)S_1+\frac{S_2+S_1^2}{2},
\end{equation*}
with $S_t = -\sum_{i=1}^4\a_i^t$ for $t=1,2$, where $\a_i \in \C$ ($1 \le i \le 4)$ is the root of $L_K(u)$.
Therefore, the average class number of $\HE_5(\F_q)$ is
\begin{equation}\label{class_num_av}
(h_K)^{\av}=\frac{1}{\# \HE_5(\F_q)}\sum_{K \in \HE_5(\F_q)} h_K = (q^2+1) + (q+1)S_1^{\av} +\frac{S_2^{\av}+(S_1^2)^{\av}}{2},
\end{equation}
where $(S_t)^{\av}$ is the average value for $S_t$'s. Since $m=1$, we have
\begin{equation}\label{s1av}
(S_1)^{\av} = (S_2)^{\av}=0
\end{equation}
by Lemma \ref{cor1}. Hence, it suffices to compute the average value on $S_1^2$. Let $J(\l_1,\l_2) = \sum_{i=1}^5a_i\z^i$, where $a_i \in \Z$, $\l_1$ (resp. $\l_2$) is a multiplicative character of $\F_q$ of order 5 (resp. 2).
By Theorem \ref{thm1}, we obtain
$$S_1^2 = (N_1-(q+1))^2 = (-a(q))^2 = \left\{\begin{array}{ll}
(\ff_{1,n})^2& \text{if $a^2-4b = \square$ in $\F_q^\times$,} \\[5pt]
(-\ff_{1,n})^2 & \text{if $a^2-4b \ne \square$ in $\F_q^\times$},
\end{array}
\right.$$
where $n$ is determined by the value of $\l_1(-4^{-1}(a^2-4b))$. Thus, we have the following:
\begin{align}\label{s12av}
(S_1^2)^{\av} &= \frac{2}{q(q-1)}\cdot \frac{q(q-1)}{2}\cdot \frac{1}{5}
\sum_{n=1}^5(\ff_{1,n})^2 = \frac{1}{5}
\sum_{n=1}^5(\ff_{1,n})^2 = \frac{1}{5}\sum_{n=1}^5 \left(5a_{5-n} - \sum_{i=1}^5a_i \right)^2,\nonumber\\[5pt]
&=\frac{1}{5}\left(20\sum_{i=1}^5a_i^2 - 10\sum_{1 \le i < j \le 5}a_ia_j \right) = 4q;
\end{align}
we use Lemma \ref{lem_5} for the last equality. Combining \eqref{class_num_av}, \eqref{s1av}, and \eqref{s12av} altogether, we get
\begin{align*}
(h_K)^{\av} = (q^2+1) + (q+1)S_1^{\av} +\frac{S_2^{\av}+(S_1^2)^{\av}}{2} =  (q^2+1) + \frac{(S_1^2)^{\av}}{2} = q^2+2q+1 = (q+1)^2.
\end{align*}
\end{proof}

\begin{thm}\label{thmD}
Let $q$ be an odd prime power and $m$ be the order of $q$ in $\F_7^\times$. The average class number of $\HE_7(\F_q) = \{k(\sqrt[7]{T^2+aT+b}) \mid a,b \in \F_q, ~a^2-4b \ne 0\}$ is given as follows.
$$(h_K)^{\av} = \left\{
\begin{array}{ll}
(q+1)^3& \text{if $m = 1, 2~~$ (i.e., $q \equiv \pm 1 \pmod 7$),}\\[5pt]
q^3+1  &\text{if $m=3, 6~~$ (i.e., $q \not\equiv \pm 1 \pmod 7$).}
\end{array}
\right.$$

Furthermore, set
$\HE_7^{\text{sq}}(\F_q) := \{k(\sqrt[7]{T^2+aT+b}) \mid a,b \in \F_q, ~a^2-4b=\square \text{ in $\F_{q^m}^\times$}\}$ and $\HE_7^{\text{non-sq}}(\F_q) := \{k(\sqrt[7]{T^2+aT+b}) \mid a,b \in \F_q, ~a^2-4b\ne\square \text{ in $\F_{q^m}$}\}$.
Let $J(\l_1,\l_2) := \sum_{i=1}^7a_i\z^i$, where $\l_1$ (resp. $\l_2$) is a multiplicative character of $\F_{q^m}$ of order 7 (resp. order 2), and $\z = e^{\frac{2\pi \imath}{7}}$.

The average class number of $\HE_7^{\text{sq}}(\F_q)$ is given as follows.
\allowdisplaybreaks
\begin{align}
(h_K)^{\av} &=\left\{
\begin{array}{ll}
\mathlarger{q^3+3q^2+2+7\bigg(2\sum_{1\le i<j<k\le 7}a_ia_ja_k-\mathfrak{A}\bigg)}& \text{if $m=1$ (i.e., $q \equiv 1 \pmod 7$),}\\[15pt]
\mathlarger{q^3+1-\frac{7a_7+1}{3}} & \text{if $m=3$ (i.e., $q \equiv 2,4 \pmod 7$).}
\end{array}
\right. \nonumber
\end{align}

The average class number of $\HE_7^{\text{non-sq}}(\F_q)$ is given as follows.
\allowdisplaybreaks
\begin{align}
(h_K)^{\av} &=\left\{
\begin{array}{ll}
\mathlarger{q^3+3q^2+6q-7\bigg(2\sum_{1\le i<j<k\le 7}a_ia_ja_k-\mathfrak{A}\bigg)} & \text{if $m=1$ (i.e., $q \equiv 1 \pmod 7$),}\\[5pt]
\mathlarger{q^3+1+\frac{7a_{7}+1}{3}} & \text{if $m=3$ (i.e., $q \equiv 2,4\pmod 7$),}
\end{array}
\right. \nonumber
\end{align}
where $\mathfrak{A}:=a_1a_2a_3 + a_1a_2a_5+a_1a_3a_5 +a_1a_3a_6+a_1a_4a_5 + a_1a_4a_6+ a_1a_4a_7+ a_1a_6a_7 + a_2a_3a_4+a_2a_3a_6+a_2a_4a_6+ a_2a_4a_7+a_2a_5a_6 + a_2a_5a_7 + a_3a_4a_5+a_3a_4a_7+a_3a_5a_7+a_3a_6a_7+a_4a_5a_6+a_5a_6a_7$.

We note that for the case where $m=2,6$, the average class numbers of $\HE^{\text{sq}}_7(\F_q)$ and $\HE^{\text{non-sq}}_7(\F_q)$ are the same by Theorem \ref{thmB}.
\end{thm}

\begin{proof}
Let $K \in \HE_7(\F_{q})$ be a hyperelliptic function field and let $L_K(u) = \sum_{i=0}^6c_iu^i \in \Z[u]$ be the $L$-polynomial of $K$. By Lemma \ref{c}, the class number of $K$ can be computed as follows:
\begin{align*}
h_K&=(q^3+1)+(q^2+1)c_1+(q+1)c_2+c_3\\
&=(q^3+1) + (q^2+1)S_1 + \frac{q+1}{2}(S_2+S_1^2) +  \frac{S_3}{3}+ \frac{S_1S_2}{2} + \frac{S_1^3}{6}
\end{align*}
with $S_t = -\sum_{i=1}^6\a_i^t$ for $1 \le t \le 3$, where $\a_i \in \C$ ($1 \le i \le 6)$ is the root of $L_K(u)$. Since the possible $m$ is 1,2,3, or 6, we only need to compute the case $m=1$ by Theorem \ref{thmB}; thus, we now assume that $m=1$. Using Lemma \ref{lem_5}, we can easily obtain that
$$(S_1^2)^{\av} = \frac{1}{7}\sum_{n=1}^7\ff_{1,n}^2 = \frac{1}{7}\left(42\sum_{i=1}^7a_i^2 - 14\sum_{1 \le i < j \le 7}a_ia_j \right) = 6q.$$
Therefore, we have
\begin{align*}
(h_K)^{\av} &= q^3+1 + \frac{q+1}{2}(S_1^2)^{\av} + \frac{(S_1S_2)^{\av}}{2} + \frac{(S_1^3)^{\av}}{6}\\
&=\left\{
\begin{array}{ll}
(q+1)^3 + \frac{1}{14}\sum_{n=1}^7(\ff_{1,n}\ff_{2,n}) +\frac{1}{42}\sum_{n=1}^7(\ff_{1,n}^3) & \text{if $K \in \HE_7^{\text{sq}}(\F_q)$,}\\[10pt]
(q+1)^3 - \frac{1}{14}\sum_{n=1}^7(\ff_{1,n}\ff_{2,n}) -\frac{1}{42}\sum_{n=1}^7(\ff_{1,n}^3) & \text{if $K \in \HE_7^{\text{non-sq}}(\F_q)$;}
\end{array}
\right.
\end{align*}
we note that $S_1^{\av} = S_2^{\av} = S_3^{\av}=0$ since $m=1$. By \eqref{lem_7}, we can easily obtain that $\ff_{1,n} = 7a_{1,\overline{7-n}}+1$. Also, by simple computation,
$$\ff_{2,n} = 7a_{2,\overline{7-2n}}-\sum_{i=1}^7a_{2,i} = 7a_{2,\overline{7-2n}}+\left(\sum_{i=1}^7a_{1,i}\right)^2 = 7a_{2,\overline{7-2n}}+1.$$
Therefore,
\begin{align}\label{comp1}
\frac{1}{14}\sum_{n=1}^7(\ff_{1,n}\ff_{2,n}) &= \frac{1}{14}\sum_{n=1}^7(7a_{1,\overline{7-n}}+1)(7a_{2,\overline{7-2n}}+1)\nonumber\\
&=\frac{7}{2}\sum_{n=1}^7(a_{1,\overline{7-n}}a_{2,\overline{7-2n}})-\frac{1}{2} = -\frac{7}{2}\sum_{n=1}^7a_{1,n}^3-\frac{1}{2}-7\mathfrak{A}.
\end{align}
Also, we have
\begin{align}\label{comp2}
\frac{1}{42}\sum_{n=1}^7(\ff_{1,n}^3) = \frac{1}{42}\sum_{n=1}^7(7a_{1,\overline{7-n}}+1)^3 = \frac{49}{6}\sum_{n=1}^7a_{1,\overline{7-n}}^3+\frac{7}{2}\sum_{n=1}^7a_{1,\overline{7-n}}^2+\frac{1}{2}\sum_{n=1}^7a_{1,\overline{7-n}}+\frac{1}{6}
\end{align}
Combining \eqref{comp1}, \eqref{comp2}, and Lemma \ref{lem_5}, we get
$$\frac{1}{14}\sum_{n=1}^7(\ff_{1,n}\ff_{2,n}) +\frac{1}{42}\sum_{n=1}^7(\ff_{1,n}^3) = -3q+1+7(2\sum_{1\le i <j<k\le 7}a_ia_ja_k-\mathfrak{A});$$
thus, the result follows.
\end{proof}

\begin{rem}
We note that from Theorems \ref{thmC} and \ref{thmD}, \mbox{the average class number of $\HE_\ell(\F_q)$ is}
\allowdisplaybreaks
\begin{align*}
(h_K)^{\av} = \left\{
\begin{array}{ll}
(q+1)^\frac{\ell-1}{2}& \text{if $q \equiv \pm 1 \pmod \ell$,}\\[5pt]
q^\frac{\ell-1}{2}+1  &\text{if $q \not \equiv \pm 1 \pmod \ell$}
\end{array}
\right.
\end{align*}
for $\ell=5$ and $7$.
We might guess that this result would be true for any odd prime $\ell$. However, in Example \ref{class_av_ex}, we find a counterexample for $\ell=11$, which is the smallest prime after the prime $7$.
\end{rem}

\vskip 0.2cm

\begin{ex}\label{class_av_ex}
Let $q=23$ and $\ell=11$. We consider a hyperelliptic curve of genus 5 and the form
$$y^{11} = x^2+ax+b \in \F_{23}[x] \text{ with $a^2-4b \ne 0$ in $\F_{23}$}.$$

We first assume that $a^2-4b$ is a square in $\F_{23}$. Depending on the value of $\l_1(-4^{-1}(a^2-4b))=\l_1(17(a^2-4b))$, we get
\begin{align*}
(S_1, S_2,S_3,S_4,S_5) \in &\{(-21,1,-219,-967,-10361), (23,-131,-175,89,-10647), \\
&(1,89,529,-2199,3081), (1,-43,-593, 2113,-2639),\\
&(23,45,-43,-1935,6073), (1,133,463,-1935,-15069),\\
&(-21,45,309,2465,7239), (-21, -87, 375, 1849, 8339)\\
&(1,-43,-329,1673,1651), (1,1,-65,-791,7921), (12, -10, -252, -362,4412)\}
\end{align*}
and
\begin{align*}
(c_1, c_2,c_3,c_4,c_5) \in &\{(
-21,221, -1627,9505, -47937), (23,199,463,-4839,-43075),\\
&(1, 45, 221, 639, 8999),( 1, -21, -219, 551, 4379),\\
&(23, 287, 2531, 17051, 91059), ( 1, 67, 221, 1915, 9065),\\
&(-21, 243, -1913, 11771, -60543),(-21, 177, -505, -2705, 27325 ),\\
&(1, -21, -131, 529, 3279), (1, 1, -21, -219, 1365), (12, 67, 144, -582, -5048)\}.
\end{align*}

Therefore, the possible class numbers are
$$\{2566663, 15380321, 7405211, 6362191, 18206639, 7703597, 2724557, 2408153, 6407203, 6713333,10667008\};$$
thus, the average class number of $\{k(\sqrt[11]{T^2+aT+b}) \mid a, b \in \F_{23}, a^2-4b = \square\}$ is $7867716.$

Using a similar computation as in the square case, the possible class numbers of $\{k(\sqrt[11]{T^2+aT+b}) \mid a, b \in \F_{23}, a^2-4b \ne \square\}$ are
$$\{16140521, 2102959, 6593269, 6025889, 2468929, 6891523, 16626787, 14642167, 5979821, 6173179,3808256\};$$
hence, the average is $7950300$. Thus, the average class number of $\HE_{11}(\F_{23})$ is
$$\frac{7867716+7950300}{2} = 7909008 \ne 24^5 = 7962624.$$
\end{ex}


\section{Implementation results}

In this section, we demonstrate some examples which use Theorem \ref{thmA} and Corollary \ref{corA}. More specifically, we compute the case where $(m,g) = (1,6), (1,20), (3,15)$ in Example \ref{ex2}, \ref{class_num_ex2}, and \ref{class_num_ex4}, respectively, where $m$ is the order of $q$ in $\F_\ell^\times$ and $g$ is the genus of the given function field. We also compute the average class numbers of hyperelliptic function fields of $\HE_5(\F_{31})$ in Example \ref{ex1} which indicates Theorem \ref{thmC}. As a final example, we calculate the average class numbers of hyperelliptic function fields of $\HE_{11}(\F_{23})$. For each example, we give detailed explanation.


In Example \ref{ex2}, we compute the trace of Frobenius of the following given hyperelliptic curve
\begin{equation}\label{ff1}
HE_{44,23}=y^{13} = x^2+44x+23 \hbox{ over $\F_{53}$,}
\end{equation}
and compute the number of $\F_{53^t}$-points (up to $t \le 13$). Using these data, we also calculate the class number of corresponding function field $\F_{53}(T, \sqrt[13]{T^2+44T+23})$.

\vskip 0.4cm

\begin{ex}\label{ex2}
We note that $x^2+44x+23 \in \F_{53}[x]$ has a square discriminant: that is, $44^2-4\times 23 = 42 = 25^2$ in $\F_{53}$. Also, the order of 53 in $\F_{13}^\times$ is 1: that is, $m=1$.

Let $\l_1$ be a multiplicative character on $\F_{53}$ of order 13 such that $\l_1(2) = \z$ with $\z = e^{\frac{2\pi \imath}{13}}$, and $\l_2$ a multiplicative character on $\F_{53}$ of order 2. Then we obtain \hbox{the following Jacobi sum}
$$J(\l_1,\l_2) = \sum_{i=1}^{13}a_{1,i}\z^i =2\z^3+2\z^4-2\z^5+2\z^6-2\z^7+2\z^8+2\z^9-2\z^{10}-4\z^{11}-1.$$
Since $\l_1(-4^{-1}(17-45)^2) = \l_1(50) = \z^4$, we take $n=4$. By Theorem \ref{thm1}, we have the following:
$$a(q) = -\ff_{1,4} = -13a_{1,13-4}+\sum_{i=1}^{13}a_{1,i} = -13\cdot 2 -1 = -27.$$
Similarly, we get
$$-J(\l_1,\l_2)^2 = \sum_{i=1}^{13}a_{2,i}\z^i =
43 + 48\z + 36\z^2+28\z^3+32\z^4+12\z^5+52\z^6+28\z^7+20\z^8+16\z^9 +20\z^{10} + 28\z^{11};$$
thus,
$$a(q^2) = -\ff_{2,4} = -13a_{2,13-2\cdot 4}+\sum_{i=1}^{13}a_{2,i} = -13\cdot 12 + 363 = 207.$$
As in Theorem \ref{thmA}, let $(-1)^{t-1}J(\l_1,\l_2)^t = \sum_{i=1}^{13}a_{t,i}\z^i$, then we have the following:
\allowdisplaybreaks
{\small{\begin{align*}
&a(q^3) = -\ff_{3,4} = -13a_{3,1}+\sum_{i=1}^{13}a_{3,i} = -13\cdot 0-261=-261 \\
& a(q^4) = -\ff_{4,4} = -13a_{4,10}+\sum_{i=1}^{13}a_{4,i} = -13 \cdot (-56) -5929 = -5201\\
& a(q^5)= -\ff_{5,4} = -13a_{5,6}+\sum_{i=1}^{13}a_{5,i} =-13\cdot (-8806)-45865 = 68613\\
&a(q^6) = -\ff_{6,4} = -13a_{6,2}+\sum_{i=1}^{13}a_{6,i} =-13 \cdot 17596+227395=-1353\\
& a(q^7) = -\ff_{7,4} = -13a_{7,11}+\sum_{i=1}^{13}a_{7,i} =-13\cdot (-125168)+608659=2235843\\
& a(q^8)= -\ff_{8,4} = -13a_{8,7}+\sum_{i=1}^{13}a_{8,i} =-13 \cdot 2763248 -8677969 = -44600193\\
&a(q^9) = -\ff_{9,4} = -13a_{9,3}+\sum_{i=1}^{13}a_{9,i} = -13 \cdot 9660658 + 119568175=-6020379\\
& a(q^{10}) = -\ff_{10,4} = -13a_{10,12}+\sum_{i=1}^{13}a_{10,i} =-13 \cdot 0-2091127013 = -2091127013\\
& a(q^{11})= -\ff_{11,4} = -13a_{11,8}+\sum_{i=1}^{13}a_{11,i} = -13 \cdot 510814814 + 3212656043 = -3427936539\\
&a(q^{12}) = -\ff_{12,4} = -13a_{12,4}+\sum_{i=1}^{13}a_{12,i} = -13 \cdot 1361625280 + 88479907463 = 70778778823\\
& a(q^{13}) = -\ff_{13,4} = -13a_{13,13}+\sum_{i=1}^{13}a_{13,i} = -13 \cdot 72795569287 +137880801031 = -808461599700.
\end{align*}}}

The number of $\F_{53^t}$-points is given in Table \ref{table1} as follows, where $1 \le t \le 13$.

{\small{
\begin{longtable}{|c|c|}
\centering
$t$ & the number of $\F_{53^t}$-points\\
\hline
1& 81\\
2& 2603\\
3& 149139\\
4& 7895683\\
5& 418126881\\
6& 22164362483\\
7& 1174708903995\\
8& 62259735011555\\
9& 3299763597822513\\
10& 174887472456640063\\
11& 9269035932800128137\\
12& 491258904185947375819\\
13& 26036721926414947795674\\
\caption{{\small{The number of $\F_{53^t}$-points of $y^{13} = x^2+44x+23$ over $\F_{53^t}$ ($1 \le t \le 13$)}}}
\label{table1}
\end{longtable}}}

Let $K=\F_{53}(T,\sqrt[13]{T^2+44T+23})$ be a hyperelliptic function field which corresponds to \eqref{ff1}. Then the genus of $K$ is $\frac{13-1}{2}=6$. We use Theorem \ref{thmB} to estimate the $L$-polynomial $L_K$ of $K$. We can easily find out that $a^2-4b = 44^2-4\times23=42 = \square$ in $\F_{53}$. Choosing $\l_1$ with $\l_1(2) = \z$, then we have $n=4$. For $1 \le t \le \frac{13-1}{2}=6$, we have the following $S_t$:
\begin{align*}
&S_1 = 27 \quad && S_2 = -207 \quad &&S_3 = 261\\
&S_4 = 5201 \quad && S_5 = -68613 \quad &&S_6 = 1353.
\end{align*}
For $1 \le i \le 6$, we compute $c_i$ as follows:
\allowdisplaybreaks
\begin{align*}
&c_0=1, \qquad c_1 = S_1 = 27, \qquad c_2 = \frac{S_2c_0+S_1c_1}{2} = \frac{-207+27^2}{2}= 261\\[3pt]
&c_3 = 573, \qquad c_4=-6577, \qquad c_5=-31251, \qquad c_6=28913.
\end{align*}
Therefore, the $L$-polynomial $L_K$ of $K$ is
$$L_K(u) = \sum_{i=0}^5c_i(u^i+53^{6-i}u^{12-i})+c_6u^6$$
and the class number $h_K$ of $K$ is
\begin{align*}
h_K &=\sum_{i=0}^5(1+53^{6-i})c_i + c_6 =35580222353.
\end{align*}
\end{ex}

\vskip 0.5cm


\begin{ex}\label{class_num_ex2}
Let $\ell=41$, $q=83$, $a=23$, and $b=13$. Then we have a function field
$$K=k(\sqrt[41]{T^2+23T+13}) ~~\text{with genus 20, where $k=\F_{83}(T)$}.$$

Since the order of 83 in $\F_{41}^\times$ is 1, we take $m=1$. Let $\l_1$ be a multiplicative character $\l_1$ of $\F_{83}$ of order 41 such that $\l_1(5) = \z$ with $\z=e^{\frac{2\pi \imath}{41}}$. Then we can compute that $\l_1(-4^{-1}(a^2-4b))=\l_1(-21\times 62) = \l_1(26)=\z^{12}$; thus, $n=12$. We now find a Jacobi sum $J(\l_1,\l_2)$, where $\l_2$ is a multiplicative character on $\F_{83}$ of order 2:
\begin{align*}
J(\l_1,\l_2) &= -2\z+2\z^2-2\z^3 -2\z^9+2\z^{13}-2\z^{15}-2\z^{16}+2\z^{18}+2\z^{19}-2\z^{20}-2\z^{24}+2\z^{27}\\
&\hskip 0.5cm +2\z^{29}+2\z^{30}+2\z^{31}+2\z^{33}-2\z^{34}-2\z^{35}+2\z^{36}-2\z^{37}-1.
\end{align*}

For $1 \le t \le \frac{41-1}{2}=20$, we have the following $S_t$:
\allowdisplaybreaks
\begin{align*}
&S_1 = -83, \qquad S_2 = -819, \qquad S_3= -5249, \qquad S_4 = -130871, \qquad S_5=-230913,\\[3pt]
&S_6=711105, \qquad S_7=67175711,\qquad S_8=-280533151,\qquad S_9=4988816449,\\[3pt]
&S_{10}=9473359141, \qquad S_{11}=277756050291, \qquad S_{12}=1667378649049,\\[3pt]
&S_{13} =12334605287727, \qquad S_{14}=-86544845385859, \qquad S_{15}=228259121069931,\\[3pt]
&S_{16}=-5492857857812351,\qquad S_{17}=93234050337263985, \qquad S_{18}=-159822977180882223,\\[3pt]
&S_{19}=-4181338526895682555, \qquad S_{20}=22581110694899544169.
\end{align*}
we note that $a^2-4b=62$ is a non-square in $\F_{83}$. For $1 \le i \le 20$, we compute $c_i$ as follows:
\allowdisplaybreaks
\begin{align*}
&c_0=1, \qquad c_1 = S_1 = -83, \qquad c_2 = \frac{S_2c_0+S_1c_1}{2} = \frac{-819+(-83)^2}{2}= 3035,\\[3pt]
&c_3= -63059, \qquad c_4 = 763257, \qquad c_5=-3400869, \qquad c_6=-64859211,\\[3pt]
&c_7=1674602523,\qquad c_8=-19630766461,\qquad c_9=117159899155,\\[3pt]
&c_{10}=188753831427, \qquad c_{11}=-11962007475163, \qquad c_{12}= 134147880700481,\\[3pt]
&c_{13} = -877906261634833, \qquad c_{14}=4818516341201719, \qquad c_{15}=-51496282549144551,\\[3pt]
&c_{16}=666040188494405001,\qquad c_{17}=-4283414065857567199, \qquad c_{18}=-27784423841699400331,\\[3pt]
&c_{19}=1020314354729137351229, \qquad c_{20}=-12533721481570442380525.
\end{align*}
Hence, the $L$-polynomial $L_K(u)$ of $K$ is
$$L_K(u) = \sum_{i=0}^{19}c_i(u^i+83^{20-i}u^{40-i})+c_{20}u^{20}$$
and the class number $h_K$ of $K$ is
\begin{align*}
h_K &=\sum_{i=0}^{19}(1+83^{20-i})c_i + c_{20} = 83141651763068478983185320840629643367.
\end{align*}
\end{ex}


\begin{ex}\label{class_num_ex4}
Let $\ell=31$, $q=5^2$, $a=\xi^{14}$, and $b=\xi^{17}$, where $\xi$ is a generator of $\F_{25}^\times$. Then we get a function field
$$K=k(\sqrt[31]{T^2+\xi^{14}T+\xi^{17}}) ~~\text{with genus 15, where $k=\F_{5^2}(T)$}.$$
Since the order of $5^2$ in $\F_{31}^\times$ is $3 < \frac{31-1}{2}$, choose a multiplicative character $\l_1$ of $\F_{5^6}$ of order 31 such that $\l_1(\eta) = \z$, where $\eta$ is a generator of $\F_{5^6}^\times$ and $\z = e^{\frac{2\pi \imath}{31}}$. Then we can compute that $\l_1(-4^{-1}(a^2-4b))= \l_1(\xi^9)= \l_1(\eta^{5859}) = \z^{31}$; thus, $n=31$. We now find a Jacobi sum $J(\l_1,\l_2)$, where $\l_2$ is a multiplicative character on $\F_{5^6}$ of order 2:
\begin{align*}
J(\l_1,\l_2) &= \sum_{i=1}^{31}a_{1,i}\z^i =20\z -16\z^2+16\z^3-44\z^4+20\z^5+36\z^6-44\z^7
-16\z^{10}+4\z^{11}-8\z^{12}\\
& \hskip 2.5cm +16\z^{13} +16\z^{15}
-24\z^{16}+8\z^{17}-24\z^{18}-16\z^{19}-44\z^{20}
-8\z^{21}+8\z^{22}\\
& \hskip 2.5cm +8\z^{23}+
4\z^{24}+20\z^{25}+
36\z^{26}+4\z^{27}-24\z^{28}-8\z^{29}+36\z^{30}+23.
\end{align*}
For $1 \le t \le 15$, it suffices to compute $S_t$ with $3 \mid t$: that is, $S_3$, $S_6$, $S_9$, $S_{12}$, and $S_{15}$. We have the following computational results:
\allowdisplaybreaks
\begin{align*}
&S_3=-31a_{1,31}+\sum_{i=1}^{31}a_{1,i}= -31 \times 23 -1 = -714,\\
&S_6=31a_{2,31}-\sum_{i=1}^{31}a_{1,i}= 31\times (-633)-50343 =-69966,\\
&S_9=-31a_{3,31}+\sum_{i=1}^{31}a_{3,i}= -31\times642747-9694197 =-29619354,\\
&S_{12}=31a_{4,31}-\sum_{i=1}^{31}a_{4,i}=31\times (- 53481073) - 61830863 = -1719744126,\\
&S_{15}=-31a_{5,31}+\sum_{i=1}^{31}a_{5,i}=-31 \times 13126395907 + (-108088958477) = -515007231594;
\end{align*}
we use the fact that $\xi^9$ is a non-square in $\F_{5^6}$.

For $1 \le i \le 15$, we compute $c_i$ as follows:
\allowdisplaybreaks
\begin{align*}
&c_3 = -238, \quad c_6= 16661, \quad c_{9}= -2762600, \quad c_{12}=511371250, \quad c_{15}=-51401812500,
\end{align*}
and $c_i=0$ for $3 \nmid i$. Finally, the $L$-polynomial $L_K$ of $K$ is
$$L_K(u)=\sum_{i=0}^4c_{3i}(u^{3i}+25^{15-3i}u^{30-3i})+c_{15}u^{15}$$
and the class number $h_K$ of $K$ is
\begin{align*}
h_K &=\sum_{i=0}^{4}(1+25^{15-3i})c_{3i} + c_{15} = 917199559306470093824.
\end{align*}
\end{ex}

\vskip 0.3cm


\begin{ex}\label{ex1}
Let $q=31$ and $\ell=5$. We consider a hyperelliptic curve of the form
$$y^5 = x^2+ax+b \in \F_{31}[x] \text{ with $a^2-4b \ne 0$ in $\F_{31}$}.$$
As in Theorem \ref{thm1}, let $\l_1$ be a multiplicative character on $\F_{31}$ of order 5 such that $\l_1(3) = \z$ with $\z = e^{\frac{2\pi \imath}{5}}$, and $\l_2$ a multiplicative character on $\F_{31}$ of order 2. \mbox{Then we can easily get that} $$J(\l_1,\l_2) = 2\z+2\z^3-4\z^4-1 \quad \text{and} \quad -J(\l_1,\l_2)^2= 16\z + 28\z^2 + 4\z^3 + 31.$$
As in Theorem \ref{thmD}, let $\HE_5(\F_{31}) = \HE_5^{\text{sq}}(\F_{31}) \sqcup \HE_5^{\text{non-sq}}(\F_{31})$. Depending on the value of $\l_1(-4^{-1}(a^2-4b))=\l_1(23(a^2-4b))$, we have
\begin{align*}
(S_1,S_2) \in  &
\left\{
\begin{array}{ll}
\{(-19,-59), ~(11,1), ~(1,-79), ~(11,61), ~(-4, 76)\} & \text{if $K \in \HE_5^{\text{sq}}(\F_{31})$,}\\[7pt]
\{(19,-59), ~(-11, 1), ~(-1, -79), ~(-11,61), ~(4,76)\} & \text{if $K \in \HE_5^{\text{non-sq}}(\F_{31})$.}
\end{array}
\right.
\end{align*}
Using the same notations as in Theorem \ref{thmA}, we get
\begin{align*}
(c_1, c_2) \in &
\left\{
\begin{array}{ll}
\{(-19,151), ~(11,91), ~(1, -39), ~(11,61), ~(-4,46)\}& \text{if $K \in \HE_5^{\text{sq}}(\F_{31})$,} \\[7pt]
\{(19,151), ~(-11, 61), ~(-1, -39), ~(-11,91), ~(4,46)\} & \text{if $K \in \HE_5^{\text{non-sq}}(\F_{31})$.}
\end{array}
\right.
\end{align*}
Therefore,
\begin{align*}
h_K \in &
\left\{
\begin{array}{ll}
\{505, 1405, 955, 1375, 880\} & \text{if $K \in \HE_5^{\text{sq}}(\F_{31})$,} \\[7pt]
\{1721, 671, 891, 701, 1136\} & \text{if $K \in \HE_5^{\text{non-sq}}(\F_{31})$;}
\end{array}
\right.
\end{align*}
thus, the average class number for a hyperelliptic function field $K$ is
\begin{align*}
h_K^{\av} = &
\left\{
\begin{array}{ll}
\frac{505+1405+955+1375+880}{5} = 1024 = (31+1)^2 & \text{if $K \in \HE_5^{\text{sq}}(\F_{31})$,} \\[7pt]
\frac{1721+671+891+701+1136}{5}= 1024 = (31+1)^2 & \text{if $K \in \HE_5^{\text{non-sq}}(\F_{31})$.}
\end{array}
\right.
\end{align*}
Further computation shows the following:
\begin{align*}
&J(\l_1,\l_2)^3=- 26\z - 72\z^2 -198\z^3 - 45, \quad &-J(\l_1,\l_2)^4 =  - 96\z- 1080\z^2 -232\z^3 - 273,&\\
&J(\l_1,\l_2)^5=- 2970\z - 1380\z^2 +2910\z^3+ 1939;&&
\end{align*}
hence,
\begin{align*}
&a(31) \in \{19, -11, -1, -11, 4\},  && a(31^3) \in \{19, -341, -211, 649, -116\},\\
&a(31^2) \in \{59,-1,79,-61,-76\},  && a(31^4) \in \{-1201, 3719, -521, -1681, -316 \}, \quad a(31^5) =-9196.
\end{align*}
Therefore, the average value on the trace of Frobenius of $HE_{a,b}$ over $\F_{31^t}$ ($1 \le t \le 5$) is
\begin{align*}
&a(31)^{\av} = a(31^2)^{\av} = a(31^3)^{\av} = a(31^4)^{\av} = 0 \text{ and } a(31^5)^{\av} = -9196.
\end{align*}
\end{ex}


\end{document}